\documentclass[11pt]{article}

\usepackage{amsmath}
\usepackage{amsthm}
\usepackage{amssymb}
\usepackage{amscd}
\usepackage{a4wide}
\usepackage{epic}
\usepackage{eepic}
\usepackage[pdftex]{graphicx}

\theoremstyle{definition}
\newtheorem{theorem}{Theorem}[section]
\newtheorem{prop}[theorem]{Proposition}
\newtheorem{lemma}[theorem]{Lemma}
\newtheorem{corollary}[theorem]{Corollary}
\newtheorem{definition}[theorem]{Definition}

\newtheorem{remark}[theorem]{Remark}

\numberwithin{equation}{section}

\newenvironment{demo}[1]{%
  \trivlist
  \item[\hskip\labelsep
        {\it #1.}]
}{%
\hfill\qedsymbol
  \endtrivlist
}

\newcommand\Nat{\mathbb{N}}
\newcommand\Int{\mathbb{Z}}
\newcommand\Rat{\mathbb{Q}}

\newcommand\Tab{\operatorname{SpTab}}
\newcommand\symp{\mathrm{sp}}
\newcommand\orth{\mathrm{o}}
\newcommand\GL{\mathbf{GL}}
\newcommand\Symp{\mathbf{Sp}}
\newcommand\Orth{\mathbf{O}}

\newcommand\AP{\mathcal{A}}
\newcommand\LL{\mathcal{L}}
\newcommand\Par{\mathcal{P}}
\newcommand\Even{\mathcal{E}}
\newcommand\Odd{\mathcal{O}}

\newcommand\Pf{\operatorname{Pf}}
\newcommand\trans{{}^t\!}

\newcommand\vectx{\boldsymbol{x}}
\newcommand\vecty{\boldsymbol{y}}
\newcommand\vectz{\boldsymbol{z}}
\newcommand\vecta{\boldsymbol{a}}
\newcommand\vectb{\boldsymbol{b}}

\newcommand\ep{\varepsilon}
\renewcommand\tilde{\widetilde}

\newcommand\wt{\operatorname{wt}}
\newcommand\pr{\operatorname{pr}}

\newcommand\TT{\mathcal{T}}

\newcommand\qbinom[2]%
{\bigg[\genfrac{}{}{0pt}{}{#1}{#2}\bigg]}

\allowdisplaybreaks

\title{
Intermediate Symplectic Characters \\
and \\
Shifted Plane Partitions of Shifted Double Staircase Shape
}

\author{
Soichi Okada%
\footnote{
Graduate School of Mathematics, Nagoya University, 
Furo-cho, Chikusa-ku, Nagoya 464-8602, Japan, 
{\tt okada@math.nagoya-u.ac.jp}
}
}

\date{}

\begin{document}

\maketitle

\begin{abstract}
We use intermediate symplectic characters to give a proof and 
variations of Hopkins' conjecture, now proved by Hopkins and 
Lai, on the number of shifted plane partitions of shifted double staircase 
shape with bounded entries.
In fact, we prove some character identities involving intermediate symplectic characters, 
and find generating functions for such shifted plane partitions.
The key ingredients of the proof are a bialternant formula for intermediate symplectic characters, 
which interpolates between those for Schur functions and symplectic characters, 
and the Ishikawa--Wakayama minor-summation formula.
\par
Mathematics Subject Classification (MSC2010): 05A15 (primary), 05E05, 05E10 (secondary)
\end{abstract}

\section{%
Introduction
}

This paper is motivated by a conjecture of Hopkins, 
now proved by Hopkins--Lai \cite{HL2021}, 
on the number of shifted plane partitions of shifted double staircase shape with bounded entries 
(see Theorem~\ref{thm:HL} below).
The goal of this paper is to give a proof and variations 
(including $q$-analogues) by using intermediate symplectic characters.

Given a strict partition $\mu = (\mu_1, \dots, \mu_l)$, 
the shifted diagram $S(\mu)$ of $\mu$ is defined to be the array of unit squares 
with $\mu_i$ squares in the $i$th row from top to bottom 
such that each row is indented by one square to the right 
with respect to the previous row.
A \emph{shifted plane partition of shape $\mu$} is a filling of the cells of $S(\mu)$ 
with nonnegative integers such that the entries in each row and column are weakly decreasing.
For a nonnegative integer $m$, we denote by $\AP^m(S(\mu))$ 
the set of all shifted plane partitions of shape $\mu$ with entries bounded by $m$.
We write $\delta_r = (r, r-1, \dots, 2, 1)$.
Hopkins and Lai \cite{HL2021} prove the following product formula 
by counting lozenge tilings of a certain region in the triangular lattice.

\begin{theorem}
\label{thm:HL}
(Hopkins--Lai \cite[Theorem~1.1]{HL2021})
For $0 \le k \le n$, the number of shifted plane partitions of shifted double staircase shape 
$\delta_n + \delta_k = (n+k, n+k-2, \dots, n-k+2, n-k, n-k-1, \dots, 2, 1)$ 
with entries bounded by $m$ is given by
\begin{equation}
\label{eq:HL}
\# \AP^m(S(\delta_n+\delta_k))
 =
\prod_{1 \le i \le j \le n}
 \frac{ m+i+j-1 }{ i+j-1 }
\prod_{1 \le i \le j \le k}
 \frac{ m+i+j }{ i+j }.
\end{equation}
\end{theorem}

In this paper, we shall prove and generalize this formula 
by establishing identities involving intermediate symplectic characters.
Our algebraic approach is inspired by the proofs 
of the $k=0$ case due to Macdonald \cite{Macdonald1995} 
and the $k=n$ case due to Proctor \cite{Proctor1983}.

If $k=0$, then shifted plane partitions of shifted staircase shape $\delta_n$ with entries bounded by $m$ 
are in one-to-one correspondence with symmetric plane partitions contained in the $n \times n \times m$ box.
Then we have two different $q$-analogues of (\ref{eq:HL});
one is the MacMahon conjecture \cite[p.~153]{MacMahon1899}, 
proved by Andrews \cite{Andrews1977a,Andrews1978} 
and Macdonald \cite{Macdonald1995}:
\begin{equation}
\label{eq:MacMahon}
\sum_{\sigma \in \AP^m(S(\delta_n))} q^{\| \sigma \|}
 =
\prod_{i=1}^n \frac{ [m/2+i-1/2] }{ [i-1/2] }
\prod_{1 \le i < j \le n} \frac{ [m+i+j-1] }{ [i+j-1] },
\end{equation}
and 
the other is the Bender--Knuth conjecture \cite[Eq.~(8)]{BK1972}, 
proved by Andrews \cite{Andrews1977b}, Gordon \cite{Gordon1983} and Macdonald \cite{Macdonald1995}:
\begin{equation}
\label{eq:BenderKnuth}
\sum_{\sigma \in \AP^m(S(\delta_n))} q^{| \sigma |}
 =
\prod_{1 \le i \le j \le n} \frac{ [m+i+j-1] }{ [i+j-1] }.
\end{equation}
Here we use the notation 
$\| \sigma \| = \frac{1}{2} \sum_{i=1}^n \sigma_{i,i} + \sum_{1 \le i < j \le n} \sigma_{i,j}$, 
$|\sigma| = \sum_{1 \le i \le j} \sigma_{i,j}$ 
and $[r] = (1-q^r)/(1-q)$.
Macdonald's proof \cite[I.5 Examples~16, 17 and 19]{Macdonald1995} of these identities 
proceeds as follows.
By transforming shifted plane partitions into semistandard tableaux, 
the left hand sides of (\ref{eq:MacMahon}) and (\ref{eq:BenderKnuth}) 
are expressed in terms of the $q$-specializations of Schur functions $s_\lambda(x_1, \dots, x_n)$:
\begin{gather*}
\sum_{\sigma \in \AP^m(S(\delta_n))} q^{\| \sigma \|}
 =
\sum_{\lambda \subset (m^n)} s_\lambda(q^{1/2}, q^{3/2}, \dots, q^{n-1/2}),
\\
\sum_{\sigma \in \AP^m(S(\delta_n))} q^{| \sigma |}
 =
\sum_{\lambda \subset (m^n)} s_\lambda(q, q^2, \dots, q^n),
\end{gather*}
where the sums are taken over all partitions $\lambda$ 
whose Young diagrams are contained in the $m \times n$ rectangle, i.e., the diagram of $(m^n)$.
In this setting, the key role is played by the following character identity:
\begin{equation}
\label{eq:Macdonald}
\sum_{\lambda \subset (m^n)} s_\lambda(x_1, \dots, x_n)
 =
\orth^B_{(m^n)}(x_1, \dots, x_n),
\end{equation}
where $\orth^B_{(m^n)}(x_1, \dots, x_n)$ is an odd orthogonal character, 
which is an irreducible character of $\tilde{\Orth}_{2n+1}$, 
the double-cover of the odd orthogonal group $\Orth_{2n+1}$.
Then we can obtain (\ref{eq:MacMahon}) and (\ref{eq:BenderKnuth}) 
by using the $q$-analogues of the Weyl dimension formula.

In a similar vein, Proctor \cite{Proctor1983} derived the case $k=n$ of (\ref{eq:HL}) 
from the character identity
\begin{equation}
\label{eq:Proctor}
\sum_{\lambda \subset (m^n)} \symp_\lambda(x_1, \dots, x_n)
 =
s_{(m^n)}(x_1, x_1^{-1}, \dots, x_n, x_n^{-1}, 1),
\end{equation}
where $\symp_\lambda(x_1, \dots, x_n)$ is a symplectic character, 
which is an irreducible character of the symplectic group $\Symp_{2n}$.

Now we explain our proof of Theorem~\ref{thm:HL}.
The main actor is a family of intermediate symplectic characters 
$\symp^{(k,n-k)}_\lambda(x_1, \dots, x_k | x_{k+1}, \dots, x_n)$, 
introduced by Proctor \cite{Proctor1988}, 
which are defined as the multivariate generating functions of $(k,n-k)$-symplectic tableaux 
(see Definition~\ref{def:int_symp}).
In the extreme cases, they reduce to the Schur functions $s_\lambda(\vectx) 
= \symp^{(0,n)}_\lambda(\vectx)$ 
and the symplectic characters $\symp_\lambda(\vectx) = \symp^{(n,0)}_\lambda(\vectx)$.
As a special case of our main theorem (Theorem~\ref{thm:main}), 
we obtain the following character identity.

\begin{theorem}
\label{thm:char_HL}
(the $a=0$ case of Theorem~\ref{thm:main} (1))
Let $0 \le k \le n$.
For a nonnegative integer $m$, we have
\begin{multline}
\label{eq:char_HL}
\sum_{\lambda \subset (m^n)}
 \symp^{(k,n-k)}_\lambda(x_1, \dots, x_k | x_{k+1}, \dots, x_n)
\\
 =
\orth^B_{((m/2)^n)}(x_1, \dots, x_n)
\cdot
\symp_{((m/2)^k)}(x_1, \dots, x_k)
\cdot
(x_{k+1} \cdots x_n)^{m/2}.
\end{multline}
\end{theorem}

Note that Equation~(\ref{eq:char_HL}) reduces to (\ref{eq:Macdonald}) and (\ref{eq:Proctor}) 
when $k=0$ and $n$ respectively (see also Corollary~\ref{cor:main_schur}). 
The proof of our main theorem (Theorem~\ref{thm:main}) is a generalization of proofs 
of (\ref{eq:Macdonald}) and (\ref{eq:Proctor}) 
provided in \cite{Okada1998}, and based on the Ishikawa--Wakayama minor summation 
formula \cite{IW1995}.
An additional key ingredient of the proof is a bialternant formula for 
intermediate symplectic characters (Theorem~\ref{thm:bialt}), 
which is another main result of this paper.

Theorem~\ref{thm:char_HL} enables us to find $q$-analogues of (\ref{eq:HL}).
Given a shifted plane partition $\sigma \in \AP^m(\delta_n+\delta_k)$, 
we define its weights $v(\sigma)$ and $w(\sigma)$ by putting
\begin{align}
\label{eq:wt1}
v(\sigma)
 &=
\left( k - \frac{1}{2} \right) t_0(\sigma)
 +
\sum_{l=0}^{n-k-1} t_l(\sigma)
 -
n t_{n-k}(\sigma)
 +
\sum_{l=n-k}^{n+k-1} (-1)^{l-n+k+1} (l-n+k) t_l(\sigma),
\\
\label{eq:wt2}
w(\sigma)
 &=
k t_0(\sigma)
 +
\sum_{l=0}^{n-k-1} t_l(\sigma)
 -
n t_{n-k}(\sigma)
 +
\sum_{l=n-k}^{n+k-1} (-1)^{l-n+k+1} (l-n+k+1) t_l(\sigma),
\end{align}
where
\begin{equation}
\label{eq:tr}
t_l(\sigma) = \sum_{i : (i,i+l) \in S(\delta_n+\delta_k)} \sigma_{i,i+l}
\end{equation}
is the $l$th trace of $\sigma$.

Since shifted plane partitions of shape $\delta_n + \delta_k$ are in bijection with 
$(k,n-k)$-symplectic tableaux (see Lemma~\ref{lem:bijection}), 
we obtain the following $q$-analogues of (\ref{eq:HL}) 
by specializing $x_i = q^i$ or $x_i = q^{i-1/2}$ in (\ref{eq:char_HL}).

\begin{corollary}
\label{cor:qHL}
For $0 \le k \le n$, 
the generating functions of shifted plane partitions of shifted double staircase shape 
$\delta_n + \delta_k$ with entries bounded by $m$ are given by
\begin{align}
\sum_{\sigma \in \AP^m(S(\delta_n+\delta_k))} q^{v(\sigma)}
 &=
\frac{1}{q^{mk^2/2}}
\prod_{i=1}^n \frac{ [m/2+i-1/2] }{ [i-1/2] }
\prod_{1 \le i < j \le n} \frac{ [m+i+j-1] }{ [i+j-1] }
\notag
\\
 &\quad\quad\times
\prod_{i=1}^k \frac{ [m/2+i] }{ [i] }
\prod_{1 \le i < j \le k} \frac{ [m+i+j] }{ [i+j] },
\label{eq:qHL1}
\\
\sum_{\sigma \in \AP^m(S(\delta_n+\delta_k))} q^{w(\sigma)}
 &=
\frac{1}{q^{mk(k+1)/2}}
\prod_{1 \le i \le j \le n} \frac{ [m+i+j-1] }{ [i+j-1] }
\prod_{1 \le i \le j \le k} \frac{ [m+i+j] }{ [i+j] }.
\label{eq:qHL2}
\end{align}
\end{corollary}

Now Theorem~\ref{thm:HL} is obtained by putting $q=1$.
Note that, when $k=0$, we have $v(\sigma) = \| \sigma \|$, $w(\sigma) = |\sigma|$, 
and (\ref{eq:qHL1}) (resp. (\ref{eq:qHL2})) becomes 
(\ref{eq:MacMahon}) (resp. (\ref{eq:BenderKnuth})).

The rest of this paper is organized as follows.
In Section~2, we recall a definition of intermediate symplectic characters 
and prove Jacobi--Trudi and bialternant formulas for them.
Our main theorem (Theorem~\ref{thm:main}, a generalization and variations 
of Theorem~\ref{thm:char_HL}) 
is stated in Section~3, and Sections~3 and 4 are devoted to its proof.
In Section~5, we apply our main theorem to find generating functions of 
shifted plane partitions of shifted double staircase shape, 
and to derive Hopkins--Lai's formula for the number of lozenge tilings of flashlight regions.
In Appendix A, we present Giambelli and dual Jacobi--Trudi formulas 
for intermediate symplectic characters.

\subsection*{%
Acknowledgements
}

This work was partially supported by 
JSPS Grants-in-Aid for Scientific Research No.~18K03208.
The author thanks S.~Hopkins for bringing the conjectured formula (\ref{eq:HL}) 
to his attention.

\section{%
Jacobi--Trudi and bialternant formulas
}

In this section, we recall a definition of intermediate symplectic characters 
and prove Jacobi--Trudi and bialternant formulas for them.

\subsection{%
Intermediate symplectic characters
}

A partition is a weakly decreasing sequence $\lambda = (\lambda_1, \lambda_2, \dots)$ 
of nonnegative integers with only finitely many nonzero entries.
A partition $\lambda$ is usually represented by its Young diagram $D(\lambda)$, 
which is the left-justified array of unit squares with $\lambda_i$ squares 
in the $i$th row.
The length of a partition $\lambda$, denoted by $l(\lambda)$, is 
the number of nonzero entries of $\lambda$.

Proctor \cite{Proctor1988} introduced the notion of intermediate symplectic tableaux 
to describe weight bases for indecomposable representations 
of the intermediate symplectic groups.
For $0 \le k \le n$, the intermediate symplectic group $\Symp_{2k,n-k}$ 
is defined to be the subgroup of the general linear group $\GL_{n+k}$ which preserves 
a skew-symmetric bilinear form of rank $2k$.
Then we have $\Symp_{0,n} = \GL_n$ and $\Symp_{2n,0} = \Symp_{2n}$. 

\begin{definition}
\label{def:int_symp}
Let $0 \le k \le n$ and $\lambda$ a partition of length $\le n$.
A \emph{$(k,n-k)$-symplectic tableau of shape $\lambda$} is a filling of the cells 
of the Young diagram $D(\lambda)$ with entries from
$$
\Gamma_{k,n-k}
=
\{ 1 < \overline{1} < 2 < \overline{2} < \dots < k < \overline{k} < k+1 < k+2 < \dots < n \}
$$
satisfying the following three conditions:
\begin{enumerate}
\item[(i)]
the entries in each row are weakly increasing;
\item[(ii)]
the entries in each column are strictly increasing;
\item[(iii)]
the entries of the $i$th row are greater than or equal to $i$.
\end{enumerate}
We denote by $\Tab^{(k,n-k)}(\lambda)$ the set of $(k,n-k)$-symplectic tableaux of shape $\lambda$.
Given a $(k,n-k)$-symplectic tableau $T$, we define
$$
\vectx^T
 =
\prod_{i=1}^k x_i^{m_T(i) - m_T(\overline{i})}
\prod_{i=k+1}^n x_i^{m_T(i)},
$$
where $\vectx = (x_1, \dots, x_n)$ are indeterminates 
and $m_T(\gamma)$ denotes the multiplicity of $\gamma \in \Gamma_{k,n-k}$ in $T$.
Then the \emph{$(k,n-k)$-symplectic character} corresponding to $\lambda$ is defined by
\begin{equation}
\label{eq:int_symp}
\symp^{(k,n-k)}_\lambda(x_1, \dots, x_k | x_{k+1}, \dots, x_n)
 =
\sum_{T \in \Tab^{(k,n-k)}(\lambda)} \vectx^T.
\end{equation}
\end{definition}

For example,
$$
T =
\begin{matrix}
\overline{1} & 2 & 3 & 3 \\
2 & \overline{2} & 4 \\
3 \\
4
\end{matrix}
$$
is a $(2,2)$-symplectic tableau of shape $(4,3,1,1)$ 
and $\vectx^T = x_1^{-1} x_2 x_3^3 x_4^2$.
The $(0,n)$-symplectic tableaux are the same as ordinary semistandard tableaux, 
while the $(n,0)$-symplectic tableaux are the same as King's symplectic tableaux \cite{King1976}.
Hence the Schur function $s_\lambda(\vectx)$ and the symplectic character $\symp_\lambda(\vectx)$ 
can be defined as the extremal cases ($k=0$ and $k=n$) of the intermediate symplectic characters:
\begin{equation}
s_\lambda(x_1, \dots, x_n)
 =
\symp^{(0,n)}_\lambda(x_1, \dots, x_n),
\quad
\symp_\lambda(x_1, \dots, x_n)
 =
\symp^{(n,0)}_\lambda(x_1, \dots, x_n).
\end{equation}
When $k=n-1$, the $(n-1,1)$-symplectic characters are called \emph{odd symplectic characters}.
The $(k,n-k)$-symplectic characters are collectively referred to as 
the \emph{intermediate symplectic characters}.

In general, the $(k,n-k)$-symplectic character $\symp^{(k,n-k)}_\lambda$ can be expressed as
$$
\symp^{(k,n-k)}_\lambda(x_1, \dots, x_k | x_{k+1}, \dots, x_n)
 =
\sum_\mu \symp_\mu(x_1, \dots, x_k) s_{\lambda/\mu}(x_{k+1}, \dots, x_n),
$$
where $\mu$ runs over all partitions of length $\le k$ such that $\mu \subset \lambda$, 
and $s_{\lambda/\mu}$ is the skew Schur function.
If $\lambda = (r^n)$ is a rectangular partition, 
then we have the following factorization formula.
\begin{prop}
\label{prop:factorization}
For a nonnegative integer $r$, we have
\begin{equation}
\label{eq:factorization}
\symp^{(k,n-k)}_{(r^n)}(x_1, \dots, x_k | x_{k+1}, \dots, x_n)
 =
\symp_{(r^k)}(x_1, \dots, x_k) \cdot (x_{k+1} \cdots x_n)^r.
\end{equation}
\end{prop}

\begin{demo}{Proof}
If $T$ is a $(k,n-k)$-symplectic tableau of shape $(r^n)$, 
then the first $k$ rows of $T$ form a $(k,0)$-symplectic tableau of shape $(r^k)$ 
and the $i$th row contains only the letter $i$ for $k+1 \le i \le n$.
Equation (\ref{eq:factorization}) follows from this observation.
\end{demo}

Schur functions and symplectic characters have several determinant representations.
(See, e.g., \cite[Lecture~24 and Appendix~A]{FH} and \cite[Appendix~2]{Proctor1994} 
for further information on classical group characters.)
The Weyl character formula can be written in the bialternant form:
\begin{gather}
\label{eq:schur_bialt}
s_\lambda(x_1, \dots, x_n)
 =
\frac{ \det \Big( x_j^{\lambda_i+n-i} \Big)_{1 \le i, j \le n} }
     { \det \Big( x_j^{n-i} \Big)_{1 \le i, j \le n} },
\\
\label{eq:symp_bialt}
\symp_\lambda(x_1, \dots, x_n)
 =
\frac{ \det \Big( x_j^{\lambda_i+n-i+1} - x_j^{-(\lambda_i+n-i+1)} \Big)_{1 \le i, j \le n} }
     { \det \Big( x_j^{n-i+1} - x_j^{-(n-i+1)} \Big)_{1 \le i, j \le n} }.
\end{gather}
Note that the denominators factor as follows:
\begin{equation}
\label{eq:schur_denom}
\det \Big( x_j^{n-i} \Big)_{1 \le i, j \le n}
 =
\prod_{1 \le i < j \le n} \big( x_i - x_j \big),
\end{equation}
and
\begin{multline}
\label{eq:symp_denom}
\det \Big( x_j^{n-i+1} - x_j^{-(n-i+1)} \Big)_{1 \le i, j \le n}
\\
 =
\prod_{i=1}^n
 \big( x_i - x_i^{-1} \big)
\prod_{1 \le i < j \le n}
 \big( x_i^{1/2} x_j^{1/2} - x_i^{-1/2} x_j^{-1/2} \big)
 \big( x_i^{1/2} x_j^{-1/2} - x_i^{-1/2} x_j^{1/2} \big).
\end{multline}
And the following determinant formulas are referred to as the Jacobi--Trudi identities:
\begin{gather}
\label{eq:schur_JT}
s_\lambda(x_1, \dots, x_n)
 =
\det \Big(
 h_{\lambda_i-i+j}(x_1, \dots, x_n)
\Big)_{1 \le i , j \le n},
\\
\label{eq:symp_JT}
\symp_\lambda(x_1, \dots, x_n)
 =
\det \left(
\begin{cases}
 h_{\lambda_i-i+1}(x_1^{\pm 1}, \dots, x_n^{\pm 1})
 &\text{if $j=1$}
\\
 h_{\lambda_i-i+j}(x_1^{\pm 1}, \dots, x_n^{\pm 1}) \\
 \quad + h_{\lambda_i-i-j+2}(x_1^{\pm 1}, \dots, x_n^{\pm 1})
 &\text{if $2 \le j \le n$}
\end{cases}
\right)_{1 \le i , j \le n},
\end{gather}
where $h_r(x_1, \dots, x_n)$ is the $r$th complete symmetric polynomial in $x_1, \dots, x_n$ 
and $h_r(x_1^{\pm 1}, \dots, x_n^{\pm 1}) \allowbreak = h_r(x_1, x_1^{-1}, \dots, x_n, x_n^{-1})$.
The aim of this section is to establish similar determinant formulas for 
intermediate symplectic characters.

\subsection{%
Jacobi--Trudi identities
}

The following Jacobi--Trudi identity 
can be obtained by using a lattice path interpretation of $(k,n-k)$-symplectic tableaux 
and the Lindstr\"om--Gessel--Viennot lemma 
(see \cite{Lindstrom1973} and \cite{GV1985,GV1989}).

\begin{prop}
\label{prop:JT1}
For a partition $\lambda$ of length $\le n$, we have
\begin{multline}
\label{eq:JT1}
\symp^{(k,n-k)}_\lambda(x_1, \dots, x_k | x_{k+1}, \dots, x_n)
\\
 =
\det \left(
 \begin{cases}
  h_{\lambda_i-i+j}(x_j^{\pm 1}, \dots, x_k^{\pm 1}, x_{k+1}, \dots, x_n)
  &\text{if $1 \le j \le k$}
  \\
  h_{\lambda_i-i+j}(x_j, \dots, x_n)
  &\text{if $k+1 \le j \le n$}
 \end{cases}
\right)_{1 \le i, j \le n},
\end{multline}
where $h_r(x_j^{\pm 1}, \dots, x_k^{\pm 1}, x_{k+1}, \dots, x_n)$ stands for 
the $r$th complete symmetric polynomial in $x_j, x_j^{-1}, \dots, x_k, x_k^{-1}, x_{k+1}, \dots, x_n$.
\end{prop}

\begin{demo}{Proof}
Let $G = (V,E)$ be the directed graph with vertex set
$$
V = \{ (i,j) \in \Int^2 : i \ge 1, \ 1 \le j \le n+k \}
$$
and edges directed from $u$ to $v$ whenever $v-u = (1,0)$ or $(0,1)$.
Given $u$, $v \in V$, we denote by $\LL(u;v)$ the set of all lattice paths from $u$ to $v$, 
i.e., all sequences $(w_0, w_1, \dots, w_r)$ of vertices of $G$ such that $w_0=u$, $w_r = v$ 
and $(w_{i-1}, w_i) \in E$ for $1 \le i \le r$.
A family $(P_1, \dots, P_n)$ of lattice paths $P_i$ is called non-intersecting 
if no two of them have a vertex in common.
We introduce an edge weight $\wt$ by putting
\begin{gather*}
\wt( (i,2j-1), (i+1,2j-1) ) = x_j,
\quad
\wt( (i,2j), (i+1,2j) ) = x_j^{-1}
\quad\text{for $1 \le j \le k$},
\\
\wt( (i,k+j), (i+1,k+j) = x_j
\quad\text{for $k+1 \le j \le n$,}
\end{gather*}
and $\wt( (i,j), (i,j+1) ) = 1$.
Then the weights of a lattice path $P = (w_0, w_1, \dots, w_r)$ and 
a family of lattice paths $(P_1, \dots, P_n)$ are defined by 
$\wt(P) = \prod_{i=1}^r \wt(w_{i-1},w_i)$ and $\wt(P_1, \dots, P_n) = \prod_{i=1}^n \wt(P_i)$ 
respectively.

For a partition $\lambda$ of length $\le n$, we put
$$
u_i
 = 
\begin{cases}
 (n-i+1,2i-1) &\text{if $1 \le i \le k$,} \\
 (n-i+1,k+i) &\text{if $k+1 \le i \le n$,}
\end{cases}
\quad
v_i = (\lambda_i+n-i+1,n+k)
\quad(1 \le i \le n)
$$
and denote by $\LL^0(u_1,\dots, u_n;v_1, \dots, v_n)$ the set of non-intersecting lattice paths 
$(P_1, \dots, P_n)$ such that $P_i \in \LL(u_i;v_i)$ for $1 \le i \le n$.
Then there is a weight-preserving bijection between $\Tab^{(k,n-k)}(\lambda)$ 
and $\LL^0(u_1, \dots, u_l; v_1, \dots, v_l)$.
See Figure~\ref{fig:1} for an example in the case where $n=4$, $k=2$ and $\lambda = (4,3,1,1)$.
\begin{figure}[htb]
$$
\begin{matrix}
\overline{1} & 2 & 3 & 3 \\
2 & \overline{2} & 4 \\
3 \\
4
\end{matrix}
\longleftrightarrow
\raisebox{-70pt}{
\includegraphics[bb=0 0 193 159]{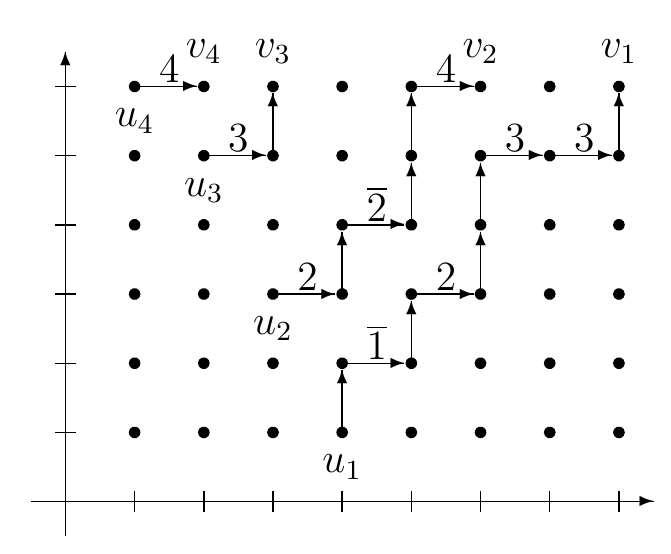}
}
$$
\caption{Lattice path interpretation}
\label{fig:1}
\end{figure}
Since the generating function of $\LL(u_j;v_i)$ is given by
$$
\sum_{P \in \LL(u_j;v_i)} \wt(P)
 =
\begin{cases}
h_{\lambda_i-i+j}(x_j^{\pm 1}, \dots, x_k^{\pm 1}, x_{k+1}, \dots, x_n)
 &\text{if $1 \le j \le k$,}
\\
h_{\lambda_i-i+j}(x_j, \dots, x_n)
 &\text{if $k+1 \le j \le n$,}
\end{cases}
$$
we can complete the proof by applying the Lindstr\"om--Gessel--Viennot lemma.
\end{demo}

By performing column operations, we can deduce the following Jacobi--Trudi-type expression 
from Proposition~\ref{prop:JT1}.

\begin{prop}
\label{prop:JT2}
For a partition $\lambda$ of length $\le n$, 
let $H^{(k,n-k)}_\lambda$ be the $n \times n$ matrix with $(i,j)$ entry given by
$$
\begin{cases}
 h_{\lambda_i-i+1}(x_1^{\pm 1}, \dots, x_k^{\pm 1}, x_{k+1}, \dots, x_n)
 &\text{if $j=1$,}
\\
 h_{\lambda_i-i+j}(x_1^{\pm 1}, \dots, x_k^{\pm 1}, x_{k+1}, \dots, x_n) \\
 \quad + h_{\lambda_i-i-j+2}(x_1^{\pm 1}, \dots, x_k^{\pm 1}, x_{k+1}, \dots, x_n)
 &\text{if $2 \le j \le k$,}
\\
 h_{\lambda_i-i+j}(x_{k+1}, \dots, x_n)
 &\text{if $k+1 \le j \le n$.}
\end{cases}
$$
Then we have
\begin{equation}
\label{eq:JT2}
\symp^{(k,n-k)}(x_1, \dots, x_k | x_{k+1}, \dots, x_n)
\\
 =
\det H^{(k,n-k)}_\lambda.
\end{equation}
\end{prop}

When $k=0$ (resp. $k=n$), this proposition reduces to the Jacobi--Trudi identity (\ref{eq:schur_JT}) 
for Schur functions (resp. (\ref{eq:symp_JT}) for symplectic characters).

\begin{demo}{Proof}
Recall that the generating function of complete symmetric functions $h_r(z_1, \dots, z_m)$ 
is given by
$$
\sum_{r \ge 0} h_r(z_1, \dots, z_m) u^r
 =
\prod_{i=1}^m \frac{ 1 }{ 1 - z_i u }.
$$
It follows from this generating function 
that
\begin{multline*}
h_r(x_{j+1}^{\pm 1}, \dots, x_k^{\pm 1}, x_{k+1}, \dots, x_n)
 + 
(x_j + x_j^{-1}) h_{r-1}(x_j^{\pm 1}, \dots, x_k^{\pm 1}, x_{k+1}, \dots, x_n)
\\
 =
h_r(x_j^{\pm 1}, \dots, x_k^{\pm 1}, x_{k+1}, \dots, x_n) 
 +
h_{r-2}(x_j^{\pm 1}, \dots, x_k^{\pm 1}, x_{k+1}, \dots, x_n)
\end{multline*}
for $1 \le j \le k$, and that
$$
h_r(x_{j+1}, \dots, x_n)
 +
x_j h_{r-1}(x_j, \dots, x_n)
 =
h_r(x_j, \dots, x_n)
$$
for $k+1 \le j \le n$.
We apply the following two types of column operations to the matrix 
on the right hand side of (\ref{eq:JT1}):
\begin{enumerate}
\item[(a)]
add the $j$th column multiplied by $x_j+x_j^{-1}$ to the $(j+1)$st column 
for $1 \le j \le k-1$;
\item[(b)]
add the $j$th column multiplied by $x_j$ to the $(j+1)$st column 
for $k+1 \le j \le n-1$.
\end{enumerate}
Then, by using the above relations, we can show
$$
\symp^{(k,n-k)}_\lambda(x_1, \dots, x_k | x_{k+1}, \dots, x_n)
 =
\det K^{(k,n-k)}_\lambda,
$$
where $K^{(k,n-k)}_\lambda$ is the $n \times n$ matrix whose $(i,j)$ entry is given by
$$
\begin{cases}
\displaystyle\sum_{l=0}^{j-1}
 \binom{j-1}{l} 
 h_{\lambda_i-i+j-2l}(x_1^{\pm 1}, \dots, x_k^{\pm 1}, x_{k+1}, \dots, x_n)
 &\text{if $1 \le j \le k$,}
\\
h_{\lambda_i-i+j}(x_{k+1}, \dots, x_n)
 &\text{if $k+1 \le j \le n$.}
\end{cases}
$$
Now, by subtracting the $j$th column multiplied by $\binom{2p+j-1}{p}$ from the $(2p+j)$th column 
for $j=1, 2, \dots, k$ and $p=1, 2, \dots, \lfloor (k-j)/2 \rfloor$, 
we conclude that 
$\det K^{(k,n-k)}_\lambda = \det H^{(k,n-k)}_\lambda$.
\end{demo}

By the same argument as in the proof of Proposition~\ref{prop:JT2}, 
we can show the following:

\begin{prop}
\label{prop:JT3}
(Proctor \cite[Proposition~3.1 and its next paragraph]{Proctor1991})
If $\lambda$ is a partition of length $l(\lambda) \le k+1$, then we have
\begin{multline}
\label{eq:JT3}
\symp^{(k,n-k)}_\lambda(x_1, \dots, x_k | x_{k+1}, \dots, x_n)
\\
 =
\det \left(
\begin{cases}
 h_{\lambda_i-i+1}(x_1^{\pm 1}, \dots, x_k^{\pm 1}, x_{k+1}, \dots, x_n)
 &\text{if $j=1$}
\\
 h_{\lambda_i-i+j}(x_1^{\pm 1}, \dots, x_k^{\pm 1}, x_{k+1}, \dots, x_n) \\
 \quad + h_{\lambda_i-i-j+2}(x_1^{\pm 1}, \dots, x_k^{\pm 1}, x_{k+1}, \dots, x_n)
 &\text{if $2 \le j \le l(\lambda)$}
\end{cases}
\right)_{1 \le i , j \le l(\lambda)}.
\end{multline}
\end{prop}

Following \cite[Definition~2.1.1]{KT1987}, 
we define the universal symplectic Schur function $s^C_\lambda(X)$ 
as a symmetric function in $X = \{ x_1, x_2, \dots \}$ by
$$
s^C_\lambda(X)
 =
\det \left(
\begin{cases}
 h_{\lambda_i-i+1}(X)
 &\text{if $j=1$}
\\
 h_{\lambda_i-i+j}(X) + h_{\lambda_i-i-j+2}(X)
 &\text{if $2 \le j \le l(\lambda)$}
\end{cases}
\right)_{1 \le i , j \le l(\lambda)},
$$
where $h_r(X)$ is the $r$th complete symmetric function in $X$.
Comparing this with (\ref{eq:JT3}), we have

\begin{corollary}
\label{cor:univ_symp}
If $l(\lambda) \le k+1$, then we have
\begin{equation}
\label{eq:univ_symp}
\symp^{(k,n-k)}_\lambda(x_1, \dots, x_k | x_{k+1}, \dots, x_n)
 =
s^C_\lambda(x_1, x_1^{-1}, \dots, x_k, x_k^{-1}, x_{k+1}, \dots, x_n, 0, 0, \dots).
\end{equation}
\end{corollary}

\subsection{%
Bialternant formulas
}

We use the Jacobi--Trudi-type identity (\ref{eq:JT2}) 
to derive bialternant formulas for $(k,n-k)$-symplectic characters, 
which interpolate between (\ref{eq:schur_bialt}) and (\ref{eq:symp_bialt}).

\begin{prop}
\label{prop:bialt}
For a partition $\lambda$ of length $\le n$, 
let $A^{(k,n-k)}_\lambda = (a_{i,j})$ be the $n \times n$ matrix with $(i,j)$ entry given by
$$
a_{i,j}
 =
\begin{cases}
 h_{\lambda_i+k-i+1}(x_j,x_{k+1}, \dots, x_n) - h_{\lambda_i+k-i+1}(x_j^{-1}, x_{k+1}, \dots, x_n)
 &\text{if $1 \le j \le k$,}
\\
 x_j^{\lambda_i+n-i}
 &\text{if $k+1 \le j \le n$.}
\end{cases}
$$
Then we have
\begin{equation}
\label{eq:bialternant1}
\symp^{(k,n-k)}_\lambda(x_1,\dots, x_k | x_{k+1}, \dots, x_n)
 =
\frac{ \det A^{(k,n-k)}_\lambda }
     { \det A^{(k,n-k)}_\emptyset },
\end{equation}
and
\begin{multline}
\label{eq:denominator1}
\det A^{(k,n-k)}_\emptyset
\\
 =
\prod_{i=1}^k
 \big( x_i - x_i^{-1} \big)
\prod_{1 \le i < j \le k}
 \big( x_i^{1/2} x_j^{1/2} - x_i^{-1/2} x_j^{-1/2} \big)
 \big( x_i^{1/2} x_j^{-1/2} - x_i^{-1/2} x_j^{1/2} \big)
\\
\times
\prod_{k+1 \le i < j \le n}
 \big( x_i - x_j \big).
\end{multline}
\end{prop}

When $k=0$ (resp. $k=n$), this theorem reduces to the bialternant formula (\ref{eq:schur_bialt}) 
(resp. (\ref{eq:symp_bialt})). 
If $k=n-1$, then we can transform this bialternant formula (\ref{eq:bialternant1}) 
into the bialternant formula obtained in \cite[Theorem~1.1]{Okada2020}.

\begin{demo}{Proof}
Let $M^{(k,n-k)} = \left( m_{i,j} \right)$ be the $n \times n$ matrices with $(i,j)$ entry given by
$$
m_{i,j}
 =
\begin{cases}
(-1)^{k-i} e_{k-i}( x_1^{\pm 1}, \dots, x_{j-1}^{\pm 1}, x_{j+1}^{\pm 1}, \dots, x_k^{\pm 1})
 \cdot (x_j - x_j^{-1})
 &\text{if $1 \le i, j \le k$,}
\\
(-1)^{n-i} e_{n-i}(x_{k+1}, \dots, x_{j-1}, x_{j+1}, \dots, x_n)
 &\text{if $k+1 \le i, j \le n$,}
\\
0
 &\text{otherwise,}
\end{cases}
$$
where $e_r(z_1, \dots, z_m)$ is the $r$th elementary symmetric polynomial in $z_1, \dots, z_m$, 
and we use the abbreviation 
$e_r(z_1^{\pm 1}, \dots, z_m^{\pm 1}) = e_r(z_1, z_1^{-1}, \dots, z_m, z_m^{-1})$.
Then we claim that
\begin{equation}
\label{eq:HM=A}
H^{(k,n-k)}_\lambda M^{(k,n-k)} = A^{(k,n-k)}_\lambda.
\end{equation}

We prove (\ref{eq:HM=A}) by computing the entries of the product 
$H^{(k,n-k)}_\lambda M^{(k,n-k)}$.
We write $\vecty = (x_1, x_1^{-1}, \dots, x_k, x_k^{-1})$ and 
$\vectz = (x_{k+1}, \dots, x_n)$.
If we put $\vecty_j = (x_1, x_1^{-1}, \dots, x_{j-1}, x_{j-1}^{-1}, x_{j+1}, \allowbreak x_{j+1}^{-1}, \dots, x_k, x_k^{-1})$, 
then we have
\begin{align*}
&
\left( \sum_{s=0}^\infty h_s(\vecty, \vectz) u^s \right)
\left( \sum_{t=0}^{2k-2} (-1)^t e_t(\vecty_j) u^t \right)
\\
&\quad=
\frac{ 1 }
     { (1 - x_j u) (1 - x_j^{-1} u) \prod_{i=k+1}^n (1 - x_i u) }
\\
&\quad=
\frac{ 1 }
     { u (x_j - x_j^{-1}) }
\left\{
 \frac{ 1 }
      { 1 - x_j u }
 -
 \frac{ 1 }
      { 1 - x_j^{-1} u }
\right\}
\frac{ 1 }
     { \prod_{i=k+1}^n (1 - x_i u) }
\\
&\quad=
\frac{ 1 }
     { u (x_j - x_j^{-1}) }
\left\{
 \sum_{s=0}^\infty h_s(x_j, \vectz) u^s
 -
 \sum_{s=0}^\infty h_s(x_j^{-1}, \vectz) u^s
\right\}.
\end{align*}
Since $e_l(\vecty_j) = e_{2k-2-l}(\vecty_j)$ for $0 \le l \le k-1$, 
we equate the coefficients of $u^{r+k-1}$ to obtain
\begin{multline}
\label{eq:HM=A1}
h_r(\vecty, \vectz) \cdot (-1)^{k-1} e_{k-1}(\vecty_j)
 +
\sum_{l=2}^k
 \left( h_{r+l-1}(\vecty,\vectz) + h_{r-l+1}(\vecty,\vectz) \right)
 \cdot (-1)^{k-l} e_{k-l}(\vecty_j)
\\
 =
\frac{ 1 }
     { x_j - x_j^{-1} }
\left( h_{r+k}(x_j, \vectz) - h_{r+k}(x_j^{-1}, \vectz) \right).
\end{multline}
Similarly, if $\vectz_j = (x_{k+1}, \dots, x_{j-1}, x_{j+1}, \dots, x_n)$, 
we have
$$
\left(
 \sum_{s=0}^\infty h_s(\vectz) u^s
\right)
\left(
 \sum_{t=0}^{n-k-1} (-1)^t e_t(\vectz_j) u^t
\right)
 =
\frac{ 1 }
     { 1 - x_j u}.
$$
Equating the coefficient of $u^{r+n}$, we obtain
\begin{equation}
\label{eq:HM=A2}
\sum_{l=k+1}^n h_{r+l}(\vectz) \cdot (-1)^{n-l} e_{n-l}(\vectz_j)
 =
x_j^{r+n}.
\end{equation}
Then the claim (\ref{eq:HM=A}) follows from (\ref{eq:HM=A1}) and (\ref{eq:HM=A2}).

Now we use (\ref{eq:HM=A}) to prove (\ref{eq:bialternant1}) and (\ref{eq:denominator1}).
Since $H^{(k,n-k)}_\emptyset$ is a upper-triangular matrix with diagonal entries $1$, 
the special case $\lambda = \emptyset$ of (\ref{eq:HM=A}) gives
\begin{equation}
\label{eq:M=A}
\det M^{(k,n-k)} = \det A^{(k,n-k)}_\emptyset.
\end{equation}
In particular, we have
$$
\det M^{(0,n)}
 = 
\det \Big( x_j^{n-i} \Big)_{1 \le i, j \le n},
\quad
\det M^{(n,0)}
 =
\det \Big( x_j^{n-i+1} - x_j^{-(n-i+1)} \Big)_{1 \le i, j \le n}.
$$
Hence, by using (\ref{eq:schur_denom}) and (\ref{eq:symp_denom}), we see that 
$\det A^{(k,n-k)}_\emptyset = \det \begin{pmatrix} M^{(k,0)} & O \\ O & M^{(0,n-k)} \end{pmatrix}$ 
is given by (\ref{eq:denominator1}).
Also it follows from (\ref{eq:JT1}), (\ref{eq:HM=A}) and (\ref{eq:M=A}) that
$$
\symp^{(k,n-k)}_\lambda(x_1, \dots, x_k|x_{k+1}, \dots, x_n)
 =
\det H^{(k,n-k)}_\lambda
 = 
\frac{ \det A^{(k,n-k)}_\lambda }
     { \det M^{(k,n-k)} }
 =
\frac{ \det A^{(k,n-k)}_\lambda }
     { \det A^{(k,n-k)}_\emptyset }.
$$
\end{demo}

For our application it is convenient to convert the bialternant formula in Proposition~\ref{prop:bialt} 
into the following form.

\begin{theorem}
\label{thm:bialt}
For a partition $\lambda$ of length $\le n$, let $\overline{A}_\lambda = (\overline{a}_{i,j})$ 
be the $n \times n$ matrix with $(i,j)$ entry given by
$$
\overline{a}_{i,j}
 =
\begin{cases}
\dfrac{ x_j^{\lambda_i+k-i+1} }
      { \prod_{l=k+1}^n (1 - x_j^{-1} x_l) }
-
\dfrac{ x_j^{-(\lambda_i+k-i+1)} }
      { \prod_{l=k+1}^n (1 - x_j x_l) }
 &\text{if $1 \le j \le k$,}
\\
x_j^{\lambda_i+n-i}
 &\text{if $k+1 \le j \le n$.}
\end{cases}
$$
Then we have
\begin{equation}
\label{eq:bialternant2}
\symp^{(k,n-k)}_\lambda(x_1,\dots, x_k | x_{k+1}, \dots, x_n)
 =
\frac{ \det \overline{A}_\lambda }
     { \det \overline{A}_\emptyset },
\end{equation}
and
\begin{multline}
\label{eq:denominator2}
\det \overline{A}^{(k,n-k)}_\emptyset
\\
 =
\prod_{i=1}^k
 \big( x_i - x_i^{-1} \big)
\prod_{1 \le i < j \le k}
 \big( x_i^{1/2} x_j^{1/2} - x_i^{-1/2} x_j^{-1/2} \big)
 \big( x_i^{1/2} x_j^{-1/2} - x_i^{-1/2} x_j^{1/2} \big)
\\
\times
\prod_{k+1 \le i < j \le n}
 \big( x_i - x_j \big).
\end{multline}
\end{theorem}

\begin{demo}{Proof}
It is enough to show that $\det \overline{A}^{(k,n-k)}_\lambda = \det A^{(k,n-k)}_\lambda$.

We write $\vectz = (x_{k+1}, \dots, x_n)$.
By the partial fraction expansion, we have
\begin{align*}
\sum_{r=0}^\infty h_{r-n+k+1}(x_j,\vectz) u^r
 &=
\frac{ u^{n-k-1} }
     { (1 - x_j u) \prod_{p=k+1}^n (1 - x_p u) }
\\
 &=
-
\sum_{p=k+1}^n
 \frac{ x_j^{-1} x_p }
      { (1 - x_j^{-1} x_p) \prod_{k+1 \le q \le n, \, q \neq p} (x_p - x_q) }
\cdot
\frac{ 1 }
     { 1 - x_p u}
\\
&\quad
+
\frac{ x_j^{-n+k+1} }
     { \prod_{q=k+1}^n (1 - x_j^{-1} x_q) }
\cdot
\frac{ 1 }
     { 1 - x_j u }.
\end{align*}
Equating the coefficients of $u^r$, we obtain
\begin{align*}
h_{r-n+k+1}(x_j,\vectz)
 &=
- \sum_{p=k+1}^n
 \frac{ x_j^{-1} }
      { (1 - x_j^{-1} x_p) \prod_{k+1 \le q \le n, \, q \neq p} (x_p - x_q) }
 x_p^{r+1}
\\
&\quad
+
\frac{ 1 }
     { \prod_{p=k+1}^n (1 - x_j^{-1} x_p) }
x_j^{r-n+k+1}.
\end{align*}
Hence we have
\begin{multline*}
h_{r-n+k+1}(x_j,\vectz) - h_{r-n+k+1}(x_j^{-1},\vectz)
\\
 =
- \sum_{p=k+1}^n
 \left\{
  \frac{ x_j^{-1} }{ (1 - x_j^{-1} x_p) }
  -
  \frac{ x_j }{ (1 - x_j x_p) }
 \right\}
 \frac{ 1 }
      { \prod_{k+1 \le q \le n, \, q \neq p} (x_p - x_q) }
 x_p^{r+1}
\\
+
\frac{ x_j^{r-n+k+1} }
     { \prod_{p=k+1}^n (1 - x_j^{-1} x_p) }
-
\frac{ x_j^{-(r-n+k+1)} }
     { \prod_{p=k+1}^n (1 - x_j x_p) }.
\end{multline*}
By using this relation, 
we can show that 
the matrix $\overline{A}^{(k,n-k)}_\lambda$ is obtained from $A^{(k,n-k)}_\lambda$ 
by adding the last $(n-k)$ columns multiplied by appropriate factors 
to the first $k$ columns.
Hence we have $\det A^{(k,n-k)}_\lambda = \det \overline{A}^{(k,n-k)}_\lambda$.
\end{demo}

\section{%
Character identities
}

In this section we state our main theorem, which gives formulas 
for certain summations of intermediate symplectic characters, 
and use the Ishikawa--Wakayama minor summation formula to express 
these summations in terms of Pfaffians.
The proof of the main theorem is completed in the next section.

\subsection{%
Main theorem
}

To state our main theorem, we introduce some notations.
For a positive integer $n$ and a nonnegative integer $m$, 
we denote by $\Par((m^n))$ the set of all partitions $\lambda$ 
whose Young diagrams are contained in the $m \times n$ rectangle, 
i.e., $l(\lambda) \le n$ and $\lambda_1 \le m$.
And we put
\begin{align*}
\Even((m^n))
 &=
\{ \lambda \in \Par((m^n)) : \text{$\lambda_i$ is even ($i=1, \dots, n$)} \},
\\
\Even'((m^n))
 &=
\{ \lambda \in \Par((m^n)) : \text{$\lambda'_i$ is even ($i=1, \dots, m$)} \},
\\
\Odd'((m^n))
 &=
\{ \lambda \in \Par((m^n)) : \text{$\lambda'_i$ is odd ($i=1, \dots, m$)} \},
\end{align*}
where $\lambda' = (\lambda'_1, \lambda'_2, \dots, \lambda'_m)$ is 
the conjugate partition of $\lambda$.

A \emph{half-partition} of length $n$ is a weakly decreasing sequence 
$\lambda = (\lambda_1, \dots, \lambda_n)$ of positive half-integers $\lambda_i \in \Nat+1/2$.
For a partition or half-partition $\lambda = (\lambda_1, \dots, \lambda_n)$, 
we define the \emph{odd orthogonal character} $\orth^B_\lambda$ 
and the \emph{even orthogonal character} $\orth^D_\lambda$ by putting
\begin{align}
\label{eq:oddorth_bialt}
\orth^B_\lambda(x_1, \dots, x_n)
 &=
\frac{ \det \Big( x_j^{\lambda_i+n-i+1/2} - x_j^{-(\lambda_i+n-i+1/2)} \Big)_{1 \le i, j \le n} }
     { \det \Big( x_j^{n-i+1/2} - x_j^{-(n-i+1/2)} \Big)_{1 \le i, j \le n} },
\\
\label{eq:evenorth_bialt}
\orth^D_\lambda(x_1, \dots, x_n)
 &=
\frac{ \det \Big( x_j^{\lambda_i+n-j} + x_j^{-(\lambda_i+n-i)} \Big)_{1 \le i, j \le n} }
     { \det \Big( x_j^{n-i} + x_j^{-(n-i)} \Big)_{1 \le i, j \le n} }
\times
\begin{cases}
 2 & \text{if $\lambda_n > 0$,} \\
 1 &\text{if $\lambda_n = 0$}
\end{cases}
\end{align}
respectively.
Note that
\begin{multline}
\label{eq:oddorth_denom}
\det \Big( x_j^{n-i+1/2} - x_j^{-(n-i+1/2)} \Big)_{1 \le i, j \le n}
\\
 =
\prod_{i=1}^n
 \big( x_i^{1/2} - x_i^{-1/2} \big)
\prod_{1 \le i < j \le n}
 \big( x_i^{1/2} x_j^{1/2} - x_i^{-1/2} x_j^{-1/2} \big)
 \big( x_i^{1/2} x_j^{-1/2} - x_i^{-1/2} x_j^{1/2} \big),
\end{multline}
and
\begin{multline}
\label{eq:evenorth_denom}
\det \Big( x_i^{n-j} + x_i^{-(n-j)} \Big)_{1 \le i, j \le n}
\\
 =
2
\prod_{1 \le i < j \le n}
 \big( x_i^{1/2} x_j^{1/2} - x_i^{-1/2} x_j^{-1/2} \big)
 \big( x_i^{1/2} x_j^{-1/2} - x_i^{-1/2} x_j^{1/2} \big).
\end{multline}
(See, for example, \cite[Appendix~A]{FH}.)
We have defined the symplectic characters $\symp_\lambda$ for partitions $\lambda$ 
by the formula (\ref{eq:symp_bialt}).
To avoid the distinction between the even and odd cases, 
we use the same formula as (\ref{eq:symp_bialt}) to define 
$\symp_\lambda(x_1, \dots, x_n)$ for a weakly decreasing sequence 
$\lambda = (\lambda_1, \dots, \lambda_n)$ of half-integers such that $\lambda_n \ge -1/2$.
Then we have
$$
\symp_{(\lambda_1, \dots, \lambda_n)}(x_1, \dots, x_n)
 =
\orth^B_{(\lambda_1+1/2,\dots,\lambda_n+1/2)}(x_1, \dots, x_n)
 \cdot 
\frac{ 1 }
     { \prod_{i=1}^n \big( x_i^{1/2} + x_i^{-1/2} \big) }.
$$
For $a \in \Int \cup (\Int + 1/2)$ 
and a partition or half-partition $\lambda = (\lambda_1, \lambda_2, \dots, \lambda_n)$, 
we write
$$
(a^n) = (\underbrace{a, \dots, a}_n),
\quad
\lambda + (a^n) = (\lambda_1+a, \lambda_2+a, \dots, \lambda_n+a).
$$

Now we can state the main theorem of this paper.

\begin{theorem}
\label{thm:main}
Let $0 \le k \le n$ and $m$ a positive integer.
\begin{enumerate}
\item[(1)]
We have
\begin{multline}
\label{eq:main1}
\sum_{\lambda \in \Par((m^n))}
 \symp^{(k,n-k)}_{\lambda+(a^n)}(x_1, \dots, x_k | x_{k+1}, \dots, x_n)
\\
 =
\orth^B_{((m/2)^n)}(x_1, \dots, x_n)
\cdot
\symp_{((m/2+a)^k)}(x_1, \dots, x_k)
\cdot
(x_{k+1} \cdots x_n)^{m/2+a}.
\end{multline}
\item[(2)]
If $m$ is even, then we have
\begin{multline}
\label{eq:main2}
\sum_{\lambda \in \Even((m^n))}
 \symp^{(k,n-k)}_{\lambda+(a^n)}(x_1, \dots, x_k | x_{k+1}, \dots, x_n)
\\
 =
\symp_{((m/2)^n)}(x_1, \dots, x_n)
\cdot
\symp_{((m/2+a)^k)}(x_1, \dots, x_k)
\cdot
(x_{k+1} \cdots x_n)^{m/2+a}.
\end{multline}
\item[(3)]
We have
\begin{multline}
\label{eq:main3}
\left(
 \sum_{\lambda \in \Even'((m^n))}
 +
 \sum_{\lambda \in \Odd'((m^n))}
\right)
 \symp^{(k,n-k)}_{\lambda+(a^n)}(x_1, \dots, x_k | x_{k+1}, \dots, x_n)
\\
 =
\orth^D_{((m/2)^n)}(x_1, \dots, x_n)
\cdot
\symp_{((m/2+a)^k)}(x_1, \dots, x_k)
\cdot
(x_{k+1} \cdots x_n)^{m/2+a}.
\end{multline}
\item[(4)]
We have
\begin{multline}
\label{eq:main4}
\left(
 \sum_{\lambda \in \Even'((m^n))}
 -
 \sum_{\lambda \in \Odd'((m^n))}
\right)
 \symp^{(k,n-k)}_{\lambda+(a^n)}(x_1, \dots, x_k | x_{k+1}, \dots, x_n)
\\
 =
(-1)^n
\cdot
\symp_{((m/2-1)^n)}(x_1, \dots, x_n)
\cdot
\orth^D_{((m/2+a+1)^k)}(x_1, \dots, x_k)
\cdot
\prod_{i=k+1}^n 
x_i^{m/2+a} (x_i-x_i^{-1}).
\end{multline}
\end{enumerate}
\end{theorem}

Note that
that (\ref{eq:main1}), (\ref{eq:main2}) and (\ref{eq:main3}) hold for $m=0$ 
by Proposition~\ref{prop:factorization} 

The extreme cases $k=0$ and $k=n$ of Theorem~\ref{thm:main} are already known 
in the literature (see \cite[Theorem~2.3 and Proof of Theorem~2.5]{Okada1998} 
and the references therein).
When $k=n-1$ or $n$ and $a=0$, we can deduce the following expressions 
of rectangular Schur functions in terms of (odd) symplectic characters, 
which also follows from \cite[Lemma~4]{Proctor1983} (see also \cite[(3.1)]{Krattenthaler1998}) 
and Corollary~\ref{cor:univ_symp}.

\begin{corollary}
\label{cor:main_schur}
For a nonnegative integer $m$, we have
\begin{gather*}
s_{(m^n)}(x_1, x_1^{-1}, \dots, x_{n-1}, x_{n-1}^{-1}, x_n,1)
 =
\sum_{\lambda \in \Par((m^n))}
 \symp^{(n-1,1)}_{\lambda}(x_1, \dots, x_{n-1} | x_n),
\\
s_{(m^n)}(x_1, x_1^{-1}, \dots, x_n, x_n^{-1}, 1)
 =
\sum_{\lambda \in \Par((m^n))}
 \symp_{\lambda}(x_1, \dots, x_n).
\end{gather*}
For a positive integer $m$, we have
\begin{gather*}
s_{(m^n)}(x_1, x_1^{-1}, \dots, x_{n-1}, x_{n-1}^{-1}, x_n)
 =
\begin{cases}
 \displaystyle\sum_{\lambda \in \Even'((m^n))}
  \symp^{(n-1,1)}_{\lambda}(x_1, \dots, x_{n-1} | x_n)
 &\text{if $n$ is even,} \\
 \displaystyle\sum_{\lambda \in \Odd'((m^n))}
  \symp^{(n-1,1)}_{\lambda}(x_1, \dots, x_{n-1} | x_n)
 &\text{if $n$ is odd,}
\end{cases}
\\
s_{(m^{n-1})}(x_1, x_1^{-1}, \dots, x_{n-1}, x_{n-1}^{-1}, x_n)
 =
\begin{cases}
 \displaystyle\sum_{\lambda \in \Odd'((m^n))}
  \symp^{(n-1,1)}_{\lambda}(x_1, \dots, x_{n-1} | x_n)
 &\text{if $n$ is even,} \\
 \displaystyle\sum_{\lambda \in \Even'((m^n))}
  \symp^{(n-1,1)}_{\lambda}(x_1, \dots, x_{n-1} | x_n)
 &\text{if $n$ is odd,}
\end{cases}
\\
s_{(m^n)}(x_1, x_1^{-1}, \dots, x_{n-1}, x_{n-1}^{-1}, x_n, x_n^{-1})
 =
\begin{cases}
 \displaystyle\sum_{\lambda \in \Even'((m^n))}
  \symp_{\lambda}(x_1, \dots, x_n)
 &\text{if $n$ is even,} \\
 \displaystyle\sum_{\lambda \in \Odd'((m^n))}
  \symp_{\lambda}(x_1, \dots, x_n)
 &\text{if $n$ is odd,}
\end{cases}
\\
s_{(m^{n-1})}(x_1, x_1^{-1}, \dots, x_{n-1}, x_{n-1}^{-1}, x_n, x_n^{-1})
 =
\begin{cases}
 \displaystyle\sum_{\lambda \in \Odd'((m^n))}
  \symp_{\lambda}(x_1, \dots, x_n)
 &\text{if $n$ is even,} \\
 \displaystyle\sum_{\lambda \in \Even'((m^n))}
  \symp_{\lambda}(x_1, \dots, x_n)
 &\text{if $n$ is odd.}
\end{cases}
\end{gather*}
\end{corollary}

\begin{demo}{Proof}
The proof is accomplished by combining Theorem~\ref{thm:main} (1), (3) and (4), 
and the following factorization formulas for rectangular Schur functions:
\begin{gather*}
s_{(m^n)}(x_1, x_1^{-1}, \dots, x_{n-1}, x_{n-1}^{-1}, x_n,1)
 =
\orth^B_{((m/2)^n)}(x_1, \dots, x_n)
\cdot
\symp_{((m/2)^{n-1})}(x_1, \dots, x_{n-1})
\cdot
x_n^{m/2},
\\
s_{(m^n)}(x_1, x_1^{-1}, \dots, x_n, x_n^{-1}, 1)
 =
\orth^B_{((m/2)^n)}(x_1, \dots, x_n)
\cdot
\symp_{((m/2)^n)}(x_1, \dots, x_n),
\end{gather*}
and, for a positive integer $m$,
\begin{align*}
&
s_{(m^n)}(x_1, x_1^{-1}, \dots, x_{n-1}, x_{n-1}^{-1}, x_n)
 +
s_{(m^{n-1})}(x_1, x_1^{-1}, \dots, x_{n-1}, x_{n-1}^{-1}, x_n)
\\
&\quad\quad
=
\orth^D_{((m/2)^n)}(x_1, \dots, x_n)
\cdot
\symp_{((m/2)^{n-1})}(x_1, \dots, x_{n-1})
\cdot
x_n^{m/2},
\\
&
s_{(m^n)}(x_1, x_1^{-1}, \dots, x_{n-1}, x_{n-1}^{-1}, x_n)
 -
s_{(m^{n-1})}(x_1, x_1^{-1}, \dots, x_{n-1}, x_{n-1}^{-1}, x_n)
\\
&\quad\quad
=
\symp_{((m/2-1)^n)}(x_1, \dots, x_n)
\cdot
\orth^D_{((m/2+1)^{n-1})}(x_1, \dots, x_{n-1})
\cdot
x_n^{m/2} (x_n-x_n^{-1}),
\\
&
s_{(m^n)}(x_1, x_1^{-1}, \dots, x_n, x_n^{-1})
 +
s_{(m^{n-1})}(x_1, x_1^{-1}, \dots, x_n, x_n^{-1})
\\
&\quad\quad
=
\orth^D_{((m/2)^n)}(x_1, \dots, x_n)
\cdot
\symp_{((m/2)^n)}(x_1, \dots, x_n),
\\
&
s_{(m^n)}(x_1, x_1^{-1}, \dots, x_n, x_n^{-1})
 -
s_{(m^{n-1})}(x_1, x_1^{-1}, \dots, x_n, x_n^{-1})
\\
&\quad\quad
=
\symp_{((m/2-1)^n)}(x_1, \dots, x_n)
\cdot
\orth_{((m/2+1)^n)}(x_1, \dots, x_n).
\end{align*}
These identities can be proved in a way similar to that of \cite{AB2019}.
\end{demo} 

The strategy of the proof of Theorem~\ref{thm:main} is the same as in \cite{Okada1998}.
Namely, 
\begin{enumerate}
\item[(a)]
First we use the bialternant formula in Theorem~\ref{thm:bialt} 
and apply the Ishikawa--Wakayama minor-summation formula (Proposition~\ref{prop:IW}) 
to express the summations in Theorem~\ref{thm:main} in terms of Pfaffians.
\item[(b)]
Next we transform the resulting Pfaffians into the products of two determinants 
(see Section~4).
\end{enumerate}

\subsection{%
Reduction to the even case
}

Now we start the proof of Theorem~\ref{thm:main}.
The following lemma enables us to reduce the proof of the odd case 
to that of the even case.

\begin{lemma}
\label{lem:reduction}
\begin{enumerate}
\item[(1)]
Let $0 < k \le n$, and $m$ a nonnegative integer.
If $\lambda = (\lambda_1, \dots, \lambda_n)$ is a partition with $\lambda_1 \le m$, 
then $x_1^m \symp^{(k,n-k)}_{(\lambda_1, \dots, \lambda_n)}(x_1, \dots, x_k | x_{k+1}, \dots, x_n)$ 
is a polynomial in $x_1$ and 
\begin{multline*}
\left[
 x_1^m \symp^{(k,n-k)}_{(\lambda_1, \dots, \lambda_n)}(x_1, \dots, x_k | x_{k+1}, \dots, x_n)
\right] \Big|_{x_1 = 0}
\\
 =
\begin{cases}
 \symp^{(k-1,n-k)}_{(\lambda_2, \dots, \lambda_n)}(x_2, \dots, x_k | x_{k+1}, \dots, x_n) 
 &\text{if $\lambda_1 = m$,} \\
 0 &\text{if $\lambda_1 < m$,}
\end{cases}
\end{multline*}
where the symbol $f|_{x_1=0}$ stands for the substitution $x_1 = 0$ in $f$.
\item[(2)]
Let $m \in \Nat \cup (\Nat + 1/2)$.
For a partition or half-partition $\lambda = (\lambda_1, \dots, \lambda_n)$ such that $\lambda_1 \le m$, 
we have
\begin{gather*}
\left[
 x_1^m \symp_{(\lambda_1, \dots, \lambda_n)}(x_1, \dots, x_n)
\right] \Big|_{x_1 = 0}
 =
\begin{cases}
 \symp_{(\lambda_2, \dots, \lambda_n)}(x_2, \dots, x_n) 
 &\text{if $\lambda_1 = m$,} \\
 0 &\text{if $\lambda_1 < m$,}
\end{cases}
\\
\left[
 x_1^m \orth^B_{(\lambda_1, \dots, \lambda_n)}(x_1, \dots, x_n)
\right] \Big|_{x_1 = 0}
 =
\begin{cases}
 \orth^B_{(\lambda_2, \dots, \lambda_n)}(x_2, \dots, x_n) 
 &\text{if $\lambda_1 = m$,} \\
 0 &\text{if $\lambda_1 < m$,}
\end{cases}
\\
\left[
 x_1^m \orth^D_{(\lambda_1, \dots, \lambda_n)}(x_1, \dots, x_n)
\right] \Big|_{x_1 = 0}
 =
\begin{cases}
 \orth^D_{(\lambda_2, \dots, \lambda_n)}(x_2, \dots, x_n) 
 &\text{if $\lambda_1 = m$,} \\
 0 &\text{if $\lambda_1 < m$.}
\end{cases}
\end{gather*}
\end{enumerate}
\end{lemma}

\begin{demo}{Proof}
We only give a proof of (1).
By the bialternant formula (\ref{eq:bialternant1}), we have
$$
x_1^m \symp^{(k,n-k)}_\lambda(x_1, \dots, x_k | x_{k+1}, \dots, x_n)
 =
\frac{ x_1^{m+k} \det A^{(k,n-k)}_\lambda }
     { x_1^k \det A^{(k,n-k)}_\emptyset }.
$$
By multiplying the first column of $A^{(k,n-k)}_\lambda = (a_{i,j})$ by $x_1^{m+k}$ 
and using $h_r(u,x_{k+1}, \dots, x_n) = \sum_{s=0}^r u^{r-s} h_s(x_{k+1}, \dots, x_n)$, 
we see that the $(i,1)$ entry becomes
$$
x_1^{m+k} a_{i,1}
 =
\sum_{s=0}^{\lambda_i+k-i+1}
 \left( x_1^{m+2k+\lambda_i-i+1-s} - x_1^{m-\lambda_i+i-1+s} \right) h_s(x_{k+1}, \dots, x_n),
$$
which is a polynomial in $x_1$.
It follows that $x_1^m \symp^{(k,n-k)}_\lambda(x_1, \dots, x_k | x_{k+1}, \dots, x_n)$ 
is a polynomial in $x_1$.
Since $\left[ x_1^{m+2k+\lambda_i-i+1 -s} - x_1^{m-\lambda_i+i-1+s} \right] \Big|_{x_1=0}
 =0$ unless $i=1$, $\lambda_1 = m$ and $s = 0$, we see that
$$
\left[
 x_1^{m+k} \det A^{(k,n-k)}_\lambda
\right]
\Big|_{x_1=0}
 =
\begin{cases}
\det \begin{pmatrix}
 -1 & * \\
 0 & A^{(k-1,n-k)}_\mu
\end{pmatrix}
 =
- \det A^{(k-1,n-k)}_\mu
 &\text{if $\lambda_1 = m$,} \\
0 &\text{if $\lambda_1 < m$,}
\end{cases}
$$
where $\mu = (\lambda_2, \dots, \lambda_n)$.
In particular, we have 
$\left[ x_1^k \det A^{(k,n-k)}_\emptyset \right] \Big|_{x_1=0}
 = - \det A^{(k-1,n-k)}_\emptyset$.
Therefore we obtain the desired identity.
\end{demo}

By using this lemma, we deduce the odd case of Theorem~\ref{thm:main} 
from the even case as follows.
Assume that Equation (\ref{eq:main1}) holds for a fixed $n$ and 
aim to prove the same equation with $n$ replaced by $n-1$.
Multiplying both sides of (\ref{eq:main1}) by $x_1^{a+m}$, we have
\begin{multline*}
\sum_{\lambda \in \Par((m^n))}
 x_1^{a+m} \symp^{(k,n-k)}_{\lambda+(a^n)}(x_1, \dots, x_k | x_{k+1}, \dots, x_n)
\\
=
x_1^{m/2} \orth^B_{((m/2)^n)}(x_1, \dots, x_n)
\cdot
x_1^{m/2+a} \symp_{((m/2+a)^k)}(x_1, \dots, x_k)
\cdot
(x_{k+1} \cdots x_n)^{m/2+a}.
\end{multline*}
Then by specializing $x_1 = 0$ and using Lemma~\ref{lem:reduction}, we obtain
\begin{multline*}
\sum_{\mu \in \Par((m^{n-1}))}
 \symp^{(k-1,n-k)}_{\mu+(a^{n-1})}(x_2, \dots, x_k | x_{k+1}, \dots, x_n)
\\
=
\orth^B_{((m/2)^{n-1})}(x_2, \dots, x_n)
\cdot
\symp_{((m/2+a)^{k-1})}(x_2, \dots, x_k)
\cdot
(x_{k+1} \cdots x_n)^{m/2+a},
\end{multline*}
which is Equation (\ref{eq:main1}) with $n$ replaced by $n-1$.
The other equations (\ref{eq:main2}), (\ref{eq:main3}) and (\ref{eq:main4}) 
are treated in the same manner.
Therefore it suffices to give a proof of Theorem~\ref{thm:main} in the case $n$ is even.

\subsection{%
Minor-summation formula
}

In the remaining of this section, we assume that $n$ is even.

Recall the Ishikawa--Wakayama minor-summation formula.
We use the notations $[n] = \{ 1, 2, \dots, n \}$ and $[0,M] = \{ 0, 1, \dots, M \}$.
Given an $n \times (M+1)$ matrix $X = (x_{ij})_{1 \le i \le n, \, 0 \le j \le M}$ 
and a subset $J = \{ j_1, \dots, j_n \}$ ($j_1 < \dots < j_n$) of column indices,  
we denote by $X([n];J)$ the $n \times n$ submatrix of $X$ 
obtained by picking up the $j_1$th, $\dots$, $j_n$th columns, 
i.e., $X([n];J) = \big( x_{i,j_q} \big)_{1 \le i, q \le n}$.
If $Y = (y_{i,j})_{0 \le i, j \le M}$ is a skew-symmetric matrix 
and $J = \{ j_1 < \dots < j_n \} \subset [0,M]$, 
then we write $Y(J) = \left( y_{j_p, j_q} \right)_{1 \le p, q \le n}$.

\begin{prop}
(Ishikawa--Wakayama \cite[Theorem~1]{IW1995})
\label{prop:IW}
Let $n$ be an even integer and $M$ a nonnegative integer.
Let $Y = (y_{ij})_{0 \le i, j \le M}$ be a skew-symmetric matrix of order $M+1$ 
and $X = (x_{ij})_{1 \le i \le n, \, 0 \le j \le M}$ an $n \times (M+1)$ matrix.
Then we have
\begin{equation}
\label{eq:IW}
\sum_{J \in \binom{[0,M]}{n}}
 \Pf Y(J) \cdot \det X \left( [n] ; J \right)
 =
\Pf \left( X Y \trans X \right),
\end{equation}
where $J$ runs over all $n$-element subsets of $[0,M]$.
\end{prop}

Given a partition $\lambda$ of length $\le n$, we put
$$
I_n(\lambda) = \{ \lambda_n, \lambda_{n-1}+1, \dots, \lambda_1+n-1 \}.
$$
Then the correspondence $\lambda \mapsto I_n(\lambda)$ gives a bijection between 
partitions $\lambda \subset (m^n)$ and $n$-element subsets of $[0,n+m-1]$.

\begin{lemma}
\label{lem:subPf}
(\cite[Lemma~3.4]{Okada1998})
Suppose that $n$ is even and $m > 0$.
Let $B = (B_{i,j})_{0 \le i, j \le n+m-1}$, 
$C = (C_{i,j})_{0 \le i, j \le n+m-1}$, 
and $D^\ep = (D^\ep_{i,j})_{0 \le i, j \le n+m-1}$, $\ep \in \{ 1, -1 \}$, 
be the $(n+m) \times (n+m)$ skew-symmetric matrices with $(i,j)$ entry 
($0 \le i < j \le n+m-1$) given by
\begin{align*}
B_{i,j}
 &= 1,
\\
C_{i,j}
 &=
\begin{cases}
 1 &\text{if $i$ is even and $j$ is odd,} \\
 0 &\text{otherwise,}
\end{cases}
\\
D^\ep_{i,j}
 &=
\begin{cases}
 1 &\text{if $j=i+1$,} \\
 \ep &\text{if $(i,j) = (0,n+m-1)$,} \\
 0 &\text{otherwise.}
\end{cases}
\end{align*}
Then, for a partition $\lambda \in \Par((m^n))$, we have
\begin{align*}
\Pf B(I_n(\lambda))
 &= 
1,
\\
\Pf C(I_n(\lambda))
 &=
\begin{cases}
 1 &\text{if $\lambda \in \Even((m^n))$,} \\
 0 &\text{otherwise,}
\end{cases}
\\
\Pf D^\ep(I_n(\lambda))
 &=
\begin{cases}
 1 &\text{if $\lambda \in \Even'((m^n))$,} \\
 \ep &\text{if $\lambda \in \Odd'((m^n))$,} \\
 0 &\text{otherwise.}
\end{cases}
\end{align*}
\end{lemma}

In order to prove Theorem~\ref{thm:main}, 
we apply the minor-summation formula (Proposition~\ref{prop:IW}) 
to the matrix $X = \left( x_{j,r} \right)_{1 \le j \le n, 0 \le r \le n+m-1}$ 
with $(j,r)$ entry given by
$$
x_{j,r}
 =
\begin{cases}
\dfrac{ x_j^{a+r-n+k+1} }
      { \prod_{l=k+1}^n (1 - x_j^{-1} x_l) }
-
\dfrac{ x_j^{-(a+r-n+k+1)} }
      { \prod_{l=k+1}^n (1 - x_j x_l) }
 &\text{if $1 \le j \le k$,}
\\
x_j^{a+r}
 &\text{if $k+1 \le j \le n$,}
\end{cases}
$$
and the skew-symmetric matrices $Y = B$, $C$ and $D^\ep$ given in Lemma~\ref{lem:subPf}.
Then it follows from the bialternant formula (\ref{eq:bialternant2}) that
$$
\symp^{(k,n-k)}_{\lambda+(a^n)}(x_1, \dots, x_k | x_{k+1}, \dots, x_n)
 =
\frac{ \det X([n];I_n(\lambda)) }
     { (-1)^{n(n-1)/2} \det \overline{A}^{(k,n-k)}_\emptyset },
$$
where $\det \overline{A}^{(k,n-k)}_\emptyset$ is given by (\ref{eq:denominator2}) 
in the factored form.
By a straightforward computation of the entries of $X Y \trans X$, 
we can see that the Pfaffian $\Pf ( X Y \trans X )$ can be expressed 
in terms of the skew symmetric matrix $Q^{n,k}(\vectx;\vecta,\vectb)$ 
introduced in the following definition.

\begin{definition}
\label{def:Q}
Let $n$ be a positive even integer and $0 \le k \le n$.
Let $\vectx = (x_1, \dots, x_n)$, $\vecta = (a_1, \dots, a_n)$ 
and $\vectb = (b_1, \dots, b_k)$ be indeterminates.
Define $Q^{n,k}(\vectx;\vecta,\vectb)$ to be the $n \times n$ skew-symmetric matrix 
with $(i,j)$ entry $Q_{i,j} = Q_{i,j}(\vectx;\vecta,\vectb)$ ($1 \le i < j \le n$) given as follows:
\begin{equation}
\label{eq:Q}
Q_{i,j}
 =
\begin{cases}
q(x_i,x_j;a_i,a_j)
\left\{
\begin{array}{r}
 \dfrac{ f(x_i^{-1}) f(x_j^{-1}) b_i b_j }{ 1 - x_i x_j }
  + \dfrac{ f(x_i^{-1}) f(x_j) b_i }{ x_j - x_i }
\\
  - \dfrac{ f(x_i) f(x_j^{-1}) b_j }{ x_j - x_i }
  - \dfrac{ f(x_i) f(x_j) }{ 1 - x_i x_j }
\end{array}
\right\}
 &\text{if $1 \le i < j \le k$,}
\\
-
q(x_i,x_j;a_i,a_j)
\left(
 \dfrac{ f(x_i^{-1}) b_i }{ 1 - x_i x_j }
 - \dfrac{ f(x_i) }{ x_j - x_i }
\right)
 &\text{if $1 \le i \le k < j \le n$,}
\\
\dfrac{ q(x_i,x_j;a_i,a_j) }
      { 1 - x_i x_j }
 &\text{if $k+1 \le i < j \le n$,}
\end{cases}
\end{equation}
where
\begin{gather}
\label{eq:q}
q(\xi, \eta ; \alpha, \beta)
 =
(\xi - \eta)(1 - \alpha \beta) + (1 - \xi \eta)(\beta - \alpha),
\\
\label{eq:f}
f(u) = \frac{ u^{n-k} }{ \prod_{l=k+1}^n (1 - u x_l) }.
\end{gather}
\end{definition}

Then we can express the summations in Theorem~\ref{thm:main} 
in terms of the Pfaffian of the skew-symmetric matrix $Q^{n,k}(\vectx;\vecta,\vectb)$.

\begin{prop}
\label{prop:sum=Pf}
Let $n$ be an even integer and $m > 0$.
We write $\vectx^r = (x_1^r, \dots, x_n^r)$ and 
$\vectx_{[k]}^r = (x_1^r, \dots, x_k^r)$.
\begin{enumerate}
\item[(1)]
We have
\begin{align*}
&
\sum_{\lambda \in \Par((m^n))}
 \symp^{(k,n-k)}_{\lambda+(a^n)}(x_1, \dots, x_k | x_{k+1}, \dots, x_n)
\\
 &=
\frac{ 1 }{ (-1)^{n(n-1)} \det \overline{A}^{(k,n-k)}_\emptyset }
\cdot
\frac{ \prod_{i=1}^k x_i^{-a} \prod_{i=k+1}^n x_i^a}
     { \prod_{i=1}^k x_i^{m+n} \prod_{i=1}^n (1 - x_i)}
\Pf Q^{n,k}(\vectx; - \vectx^{m+n}, - \vectx_{[k]}^{2a+m+n+1}).
\end{align*}
\item[(2)]
If $m$ is even, then we have
\begin{align*}
&
\sum_{\lambda \in \Even((m^n))}
 \symp^{(k,n-k)}_{\lambda+(a^n)}(x_1, \dots, x_k | x_{k+1}, \dots, x_n)
\\
 &=
\frac{ 1 }{ (-1)^{n(n-1)} \det \overline{A}^{(k,n-k)}_\emptyset }
\cdot
\frac{ \prod_{i=1}^k x_i^{-a} \prod_{i=k+1}^n x_i^a}
     { \prod_{i=1}^k x_i^{m+n} \prod_{i=1}^n (1 - x_i^2)}
\Pf Q^{n,k}(\vectx; - \vectx^{m+n+1}, - \vectx_{[k]}^{2a+m+n+1}).
\end{align*}
\item[(3)]
We have
\begin{align*}
&
\left(
 \sum_{\lambda \in \Even'((m^n))}
 +
 \sum_{\lambda \in \Odd'((m^n))}
\right)
 \symp^{(k,n-k)}_{\lambda+(a^n)}(x_1, \dots, x_k | x_{k+1}, \dots, x_n)
\\
 &=
\frac{ 1 }{ (-1)^{n(n-1)} \det \overline{A}^{(k,n-k)}_\emptyset }
\cdot
\frac{ \prod_{i=1}^k x_i^{-a} \prod_{i=k+1}^n x_i^a }
     { \prod_{i=1}^k x_i^{m+n} }
\Pf Q^{n,k}(\vectx; \vectx^{m+n-1}, - \vectx_{[k]}^{2a+m+n+1}).
\end{align*}
\item[(4)]
We have
\begin{align*}
&
\left(
 \sum_{\lambda \in \Even'((m^n))}
 -
 \sum_{\lambda \in \Odd'((m^n))}
\right)
 \symp^{(k,n-k)}_{\lambda+(a^n)}(x_1, \dots, x_k | x_{k+1}, \dots, x_n)
\\
 &=
\frac{ 1 }{ (-1)^{n(n-1)} \det \overline{A}^{(k,n-k)}_\emptyset }
\cdot
\frac{ \prod_{i=1}^k x_i^{-a} \prod_{i=k+1}^n x_i^a }
     { \prod_{i=1}^k x_i^{m+n} }
\Pf Q^{n,k}(\vectx; - \vectx^{m+n-1}, \vectx_{[k]}^{2a+m+n+1}).
\end{align*}
\end{enumerate}
\end{prop}

\begin{demo}{Proof}
As the proofs are similar, we give a sketch of the proof of (1).
By using Theorem~\ref{thm:bialt}, Lemma~\ref{lem:subPf} and Proposition~\ref{prop:IW}, 
we obtain
\begin{align*}
\sum_{\lambda \in \Par((m^n))}
 \symp^{(k,n-k)}_{\lambda+(a^n)}(x_1, \dots, x_k | x_{k+1}, \dots, x_n)
 &=
\sum_{ I \in \binom{[0,n+m-1]}{n} }
 \Pf B(I)
 \frac{ \det X([n];I) }
      { (-1)^{n(n-1)/2} \det \overline{A}^{(k,n-k)}_\emptyset }
\\
 &=
\frac{ 1 }
     { (-1)^{n(n-1)/2} \det \overline{A}^{(k,n-k)}_\emptyset }
\cdot
\Pf \left( X B \trans X \right).
\end{align*}
A direct computation gives us
$$
\sum_{0 \le r < s \le m+n-1} B_{r,s} \left( x^r y^s - x^s y^r \right)
 =
\frac{ q(x,y;-x^{m+n},-y^{m+n}) }
     { (1-x) (1-y) (1 - xy) }.
$$
By replacing $x$ with $x^{-1}$ or/and $y$ with $y^{-1}$, we obtain
\begin{align*}
\sum_{0 \le r < s \le m+n-1} B_{r,s} \left( x^{-r} y^s - x^{-s} y^r \right)
 &=
\frac{ q(x,y;-x^{m+n},-y^{m+n}) }
     { x^{m+n-1} (1-x) (1-y) (y-x) },
\\
\sum_{0 \le r < s \le m+n-1} B_{r,s} \left( x^r y^{-s} - x^s y^{-r} \right)
 &=
- 
\frac{ q(x,y;-x^{m+n},-y^{m+n}) }
     { y^{m+n-1} (1-x) (1-y) (y-x) },
\\
\sum_{0 \le r < s \le m+n-1} B_{r,s} \left( x^{-r} y^{-s} - x^{-s} y^{-r} \right)
 &=
-
\frac{ q(x,y;-x^{m+n},-y^{m+n}) }
     { x^{m+n-1} y^{m+n-1} (1-x) (1-y) (1 - xy) }.
\end{align*}
Using these relations, we can explicitly compute the entries of $X B \trans X$.
By a straightforward computation, we see that the $(i,j)$ entry of $X B \trans X$ 
is equal to
$$
Q_{i,j}(\vectx;-\vectx^{m+n},-\vectx_{[k]}^{2a+m+n+1})
\times
\begin{cases}
 \dfrac{ x_i^{-a} x_j^{-a} }
       { x_i^{m+n} x_j^{m+n} (1 - x_i) (1 - x_j) }
 &\text{if $1 \le i < j \le k$,}
\\
 \dfrac{ x_i^{-a} x_j^a }
       { x_i^{m+n} (1 - x_i) (1 - x_j) }
 &\text{if $1 \le i \le k < j \le n$,}
\\
 \dfrac{ x_i^a x_j^a }
       { (1 - x_i) (1 - x_j) }
 &\text{if $k+1 \le i < j \le n$.}
\end{cases}
$$
Hence, by using the multilinearilty of Pfaffians, we obtain
$$
\Pf \left( X B \trans X \right)
 =
\frac{ \prod_{i=1}^k x_i^{-a} \prod_{i=k+1}^n x_i^a }
     { \prod_{i=1}^k x_i^{m+n} \prod_{i=1}^n (1 - x_i) }
\Pf Q^{n,k}(\vectx;-\vectx^{m+n}, - \vectx^{2a+m+n+1}).
$$
\end{demo}

\section{%
Pfaffian identity
}

In this section, we complete the proof of Theorem~\ref{thm:main} 
by establishing a Pfaffian identity.

\subsection{
Pfaffian identity
}

By Proposition~\ref{prop:sum=Pf}, 
we need to evaluate the Pfaffian $\Pf Q^{n,k}(\vectx;\vecta,\vectb)$ 
in order to prove Theorem~\ref{thm:main}.

\begin{definition}
\label{def:WU}
For indeterminates $\vectx = (x_1, \dots, x_n)$ and $\vecta = (a_1, \dots, a_n)$, 
let $W^n(\vectx;\vecta)$ be the $n \times n$ matrix with the $i$th row
$$
\begin{pmatrix}
1 + a_i x_i^{n-1} & x_i + a_i x_i^{n-2} & \dots & x_i^{n-1} + a_i
\end{pmatrix}.
$$
For indeterminates $\vecty = (y_1, \dots, y_k)$ and $\vectb = (b_1, \dots, b_k)$, 
let $U^{k,n}(\vecty;\vectb)$ be the $k \times k$ matrix with the $i$th row
$$
\begin{pmatrix}
 y_i^{n-k} + b_i y_i^{k-1} & y_i^{n-k+1} + b_i y_i^{k-2} & \dots & y_i^{n-1} + b_i
\end{pmatrix}.
$$
\end{definition}

Then the Pfaffian $\Pf Q^{n,k}(\vectx;\vecta,\vectb)$ is evaluated as follows.

\begin{theorem}
\label{thm:Pf=det*det}
Let $n$ be an even integer and $0 \le k \le n$.
Let $\vectx = (x_1, \dots, x_n)$, $\vecta = (a_1, \dots, a_n)$ 
and $\vectb = (b_1, \dots, b_k)$ be indeterminates.
Then we have
\begin{multline}
\label{eq:Pf=det*det}
\Pf Q^{n,k}(\vectx;\vecta,\vectb)
\\
 =
\frac{ (-1)^{k(k-1)/2}
       \det W^n(\vectx;\vecta) \cdot \det U^{k,n}(\vectx_{[k]};\vectb) }
     { \prod_{1 \le i < j \le k} (x_j - x_i)(1 - x_j x_j)
       \prod_{i=1}^k \prod_{j=k+1}^n (x_j - x_i)(1 - x_i x_j)
       \prod_{k+1 \le i < j \le n} (1 - x_i x_j) },
\end{multline}
where $\vectx_{[k]} = (x_1, \dots, x_k)$.
\end{theorem}

In the extreme cases $k=0$ and $k=n$, the skew-symmetric matrix $Q^{n,k}(\vectx;\vecta,\vectb)$ 
becomes
$$
Q^{n,0}(\vectx;\vecta)
 =
\left(
 \frac{ q(x_i,x_j;a_i,a_j) }
      { 1 - x_i x_j }
\right)_{1 \le i, j \le n}
$$
and
$$
Q^{n,n}(\vectx;\vecta,\vectb)
 =
\left(
-
\frac{ q(x_i,x_j;a_i,a_j) q(x_i,x_j;b_i,b_j) }
     { (x_j - x_i)(1 - x_i x_j) }
\right),
$$
respectively.
In these cases, the assertions of Theorem~\ref{thm:Pf=det*det} are established 
in \cite[Corollary~4.6 and Theorem~4.4]{Okada1998}.

We postpone the proof of this theorem to the next subsection 
and first finish the proof of Theorem~\ref{thm:main}.

\begin{demo}{Proof of Theorem~\ref{thm:main}}
By the argument in Subsection~3.2, 
we may assume that $n$ is even.

By using the bialternant formulas (\ref{eq:symp_bialt}), (\ref{eq:oddorth_bialt}), 
(\ref{eq:evenorth_bialt}), and the denominator formulas 
(\ref{eq:symp_denom}), (\ref{eq:oddorth_denom}), (\ref{eq:evenorth_denom}), 
we can see
\begin{align*}
\det W^n(\vectx;-\vectx^{m+n})
 &=
\prod_{i=1}^n x_i^{m/2} \prod_{i=1}^n (1-x_i) \prod_{1 \le i < j \le n} (x_j-x_i)(1-x_ix_j)
\cdot
\orth^B_{((m/2)^n)}(\vectx),
\\
\det W^n(\vectx;-\vectx^{m+n+1})
 &=
\prod_{i=1}^n x_i^{m/2} \prod_{i=1}^n (1-x_i^2) \prod_{1 \le i < j \le n} (x_j-x_i)(1-x_ix_j)
\cdot
\symp^B_{((m/2)^n)}(\vectx),
\\
\det W^n(\vectx;\vectx^{m+n-1})
 &=
\prod_{i=1}^n x_i^{m/2} \prod_{1 \le i < j \le n} (x_j-x_i)(1-x_ix_j)
\cdot
\orth^D_{((m/2)^n)}(\vectx),
\end{align*}
and
\begin{align*}
&
\det U^{k,n}(\vectx;-\vectx^{m+2a+n+1})
\\
 &\quad
=
\prod_{i=1}^k x_i^{m/2+a+n-k}
\prod_{i=1}^k (1-x_i^2) \prod_{1 \le i < j \le k} (x_j-x_i)(1-x_ix_j)
\cdot
\symp_{((m/2+a)^k)}(x_1, \dots, x_k),
\\
&
\det U^{k,n}(\vectx;\vectx^{m+2a+n+1})
\\
 &\quad
=
\prod_{i=1}^k x_i^{m/2+a+n-k+1}
\prod_{1 \le i < j \le k} (x_j-x_i)(1-x_ix_j)
\cdot
\orth^D_{((m/2+a+1)^k)}(x_1, \dots, x_k).
\end{align*}
Combining these relations together with Proposition~\ref{prop:sum=Pf} and Theorem~\ref{thm:Pf=det*det}, 
we arrive at the desired identities.
\end{demo}

\subsection{%
Proof of Theorem~\ref{thm:Pf=det*det}
}

It remains to prove Theorem~\ref{thm:Pf=det*det}.

Since both sides of (\ref{eq:Pf=det*det}) have degree at most one in each of the variables 
$a_1, \dots, a_n, b_1, \dots, b_k$, 
it is enough to show that the coefficients of $\vecta^I \vectb^J = \prod_{i \in I} a_i \prod_{j \in J} b_j$ 
are the same on both sides for any subsets $I \subset [n]$ and $J \subset [k]$.
We denote by $L(I,J)$ and $R(I,J)$ the corresponding coefficients on 
the left and right hand sides respectively.

First we compute the coefficient $R(I,J)$.
We put
\begin{equation}
\label{eq:D}
\begin{aligned}
D^+_n(I)
 &= 
\{ (i,j) \in [n] \times [n] : i < j \}
\cap
\big[ \big( I \times I \big) \cup \big( ([n] \setminus I) \times ([n] \setminus I) \big) \big],
\\
D^-_n(I)
 &=
\{ (i,j) \in [n] \times [n] : i < j \}
\cap
\big[ \big( I \times ([n] \setminus I) \big) \cup \big( ([n] \setminus I) \times I \big) \big],
\\
D^+_k(J)
 &= 
\{ (i,j) \in [k] \times [k] : i < j \}
\cap
\big[ \big( J \times J \big) \cup \big( ([k] \setminus J) \times ([k] \setminus J) \big) \big],
\\
D^-_k(J)
 &=
\{ (i,j) \in [k] \times [k] : i < j \}
\cap
\big[ \big( J \times ([k] \setminus J) \big) \cup \big( ([k] \setminus J) \times J \big) \big].
\end{aligned}
\end{equation}
Then we can see that the coefficient of $\vecta^I$ in $\det W^n(\vectx;\vecta)$ is equal to
\begin{equation}
\label{eq:coeff_W}
\det \left(
 \begin{cases}
  x_i^{n-j} &\text{if $i \in I$} \\
  x_i^{j-1} &\text{if $i \not\in I$}
 \end{cases}
\right)_{1 \le i, j \le n}
 =
(-1)^{n \# I - \Sigma(I)}
\prod_{(i,j) \in D^+_n(I)} (x_j - x_i)
\prod_{(i,j) \in D^-_n(I)} (1 - x_i x_j),
\end{equation}
where $\Sigma(I) = \sum_{i \in I} i$, 
and that the coefficient of $\vectb^J$ in $\det U^{k,n}(\vectx_{[k]};\vectb)$ is equal to
\begin{multline}
\label{eq:coeff_U}
\det \left(
 \begin{cases}
  x_i^{k-j} &\text{if $i \in J$} \\
  x_i^{n-k+j-1} &\text{if $i \not\in J$}
 \end{cases}
\right)_{1 \le i, j \le k}
\\
 =
\prod_{i \in [k] \setminus J} x_i^{n-k}
\cdot
(-1)^{k \# J - \Sigma(J)}
\prod_{(i,j) \in D^+_k(J)} (x_j - x_i)
\prod_{(i,j) \in D^-_k(J)} (1 - x_i x_j),
\end{multline}
where $\Sigma(J) = \sum_{j \in J} j$.
We put
\begin{equation}
\label{eq:T}
\begin{aligned}
T^{(1)}_{n,k}
 &=
\{ (i,j) : 1 \le i < j \le k \},
\\
T^{(2)}_{n,k}
 &=
\{ (i,j) : 1 \le i \le k, \, k+1 \le j \le n \},
\\
T^{(3)}_{n,k}
 &=
\{ (i,j) : k+1 \le i < j \le n \}.
\end{aligned}
\end{equation}
Then the denominator on the right hand side of (\ref{eq:Pf=det*det}) is written as
$$
\prod_{(i,j) \in T^{(1)}_{n,k}} (x_j - x_i)(1 - x_j x_j)
\prod_{(i,j) \in T^{(2)}_{n,k}} (x_j - x_i)(1 - x_i x_j)
\prod_{(i,j) \in T^{(3)}_{n,k}} (1 - x_i x_j).
$$
Since we have
$$
D^+_n(I) \sqcup D^-_n(I) = T^{(1)}_{n,k} \sqcup T^{(2)}_{n,k} \sqcup T^{(3)}_{n,k},
\quad
D^+_k(J) \sqcup D^-_k(J) = T^{(1)}_{n,k},
$$
it follows from (\ref{eq:coeff_W}) and (\ref{eq:coeff_U}) 
that the coefficient of $\vecta^I \vectb^J$ in the right hand side of (\ref{eq:Pf=det*det}) is given by
\begin{align}
\label{eq:R}
R(I,J)
 &=
(-1)^{k(k-1)/2 + n \# I + \Sigma(I) + k \# J + \Sigma(J)}
\prod_{i \in [k] \setminus J} x_i^{n-k}
\notag
\\
 &\quad\times
\prod_{(i,j) \in D^+_n(I) \cap D^+_k(J)} \frac{ x_j - x_i }{ 1 - x_i x_j }
\prod_{(i,j) \in D^-_n(I) \cap D^-_k(J)} \frac{ 1 - x_i x_j }{ x_j - x_i }
\notag
\\
 &\quad\times
\prod_{(i,j) \in T_{n,k}^{(2)} \cap D^+_n(I)} \frac{ 1 }{1 - x_i x_j}
\prod_{(i,j) \in T_{n,k}^{(2)} \cap D^-_n(I)} \frac{ 1 }{x_j - x_i}
\notag
\\
 &\quad\times
\prod_{(i,j) \in T_{n,k}^{(3)} \cap D^+_n(I)} \frac{x_j - x_i}{1 - x_i x_j}.
\end{align}

Next we compute the coefficient $L(I,J)$ on the left hand side of (\ref{eq:Pf=det*det}).
By the multilinearity of Pfaffians, we see that $L(I,J)$ is equal 
to the Pfaffian of the skew-symmetric matrix $X(I,J)$, 
whose $(i,j)$ entry $X_{i,j}$, $i<j$, is given as follows:
\begin{enumerate}
\item[(1)]
If $1 \le i < j \le k$, then
$$
X_{i,j}
 =
\left\{
 \begin{array}{ll}
  -(x_j - x_i) & \text{if $i \in I$ and $j \in I$} \\
  -(1 - x_i x_j) & \text{if $i \in I$ and $j \not\in I$} \\
  1 - x_i x_j & \text{if $i \not\in I$ and $j \in I$} \\
  x_j - x_i & \text{if $i \not\in I$ and $j \not\in I$}
 \end{array}
\right\}
\times
\left\{
 \begin{array}{ll}
  \dfrac{ f(x_i^{-1}) f(x_j^{-1}) }{ 1 - x_i x_j } &\text{if $i \in J$ and $j \in J$} \\
  \dfrac{ f(x_i^{-1}) f(x_j) }{ x_j - x_i } &\text{if $i \in J$ and $j \not\in J$} \\
  - \dfrac{ f(x_i) f(x_j^{-1}) }{ x_j - x_i } &\text{if $i \not\in J$ and $j \in J$} \\
  - \dfrac{ f(x_i) f(x_j) }{ 1 - x_i x_j } &\text{if $i \not\in J$ and $j \not\in J$} \\
 \end{array}
\right\}.
$$
\item[(2)]
If $1 \le i \le k < j \le n$, then
$$
X_{i,j}
 =
\left\{
 \begin{array}{ll}
  -(x_j - x_i) & \text{if $i \in I$ and $j \in I$} \\
  -(1 - x_i x_j) & \text{if $i \in I$ and $j \not\in I$} \\
  1 - x_i x_j & \text{if $i \not\in I$ and $j \in I$} \\
  x_j - x_i & \text{if $i \not\in I$ and $j \not\in I$}
 \end{array}
\right\}
\times
\left\{
 \begin{array}{ll}
  - \dfrac{ f(x_i^{-1}) }{ 1 - x_i x_j } &\text{if $i \in J$} \\
  \dfrac{ f(x_i) }{ x_j - x_i } &\text{if $i \not\in J$}
 \end{array}
\right\}.
$$
\item[(3)]
If $k+1 \le i < j \le n$, then
$$
X_{i,j}
 =
\left\{
 \begin{array}{ll}
  - \dfrac{ x_j - x_i }{ 1 - x_i x_j } & \text{if $i \in I$ and $j \in I$} \\
  - 1 & \text{if $i \in I$ and $j \not\in I$} \\
  1 & \text{if $i \not\in I$ and $j \in I$} \\
  \dfrac{ x_j - x_i }{ 1 - x_i x_j } & \text{if $i \not\in I$ and $j \not\in I$}
 \end{array}
\right.
$$
\end{enumerate}
We put
\begin{equation}
\label{eq:C}
\begin{alignedat}{2}
C_1 &= \{ i : 1 \le i \le k, \, i \in I, \, i \in J \},
&\quad
C_2 &= \{ i : 1 \le i \le k, \, i \in I, \, i \not\in J \},
\\
C_3 &= \{ i : 1 \le i \le k, \, i \not\in I, \, i \in J \},
&\quad
C_4 &= \{ i : 1 \le i \le k, \, i \not\in I, \, i \not\in J \},
\\
C_5 &= \{ i : k+1 \le i \le n, \, i \in I \},
&\quad
C_6 &= \{ i : k+1 \le i \le n, \, i \not\in I \}.
\end{alignedat}
\end{equation}
By pulling out the factor $f(x_i^{-1})$ from the $i$th row/column with $i \in J$ 
and $f(x_i)$ from the $i$th row/column with $i \in [k] \setminus J$, 
and then by multiplying the last $(n-k)$ rows/columns by $-1$, 
we obtain
\begin{equation}
\label{eq:L1}
L(I,J)
 =
\Pf X(I,J)
 =
\prod_{i \in J} f(x_i^{-1}) \prod_{i \in [k] \setminus J} f(x_i)
\cdot
(-1)^{n-k}
\cdot
\Pf X'(I,J),
\end{equation}
where the entries $X'_{i,j}$ of the skew-symmetric matrix $X'(I,J)$ are given by
$$
X'_{i,j}
 =
\begin{cases}
 - \dfrac{ x_j - x_i }{ 1 - x_i x_j }
 &\text{if $(i,j) \in (C_1 \times C_1) \cup (C_1 \times C_5) \cup (C_4 \times C_4) \cup (C_5 \times C_5)$,}
 \\
 \dfrac{ x_j - x_i }{ 1 - x_i x_j }
 &\text{if $(i,j) \in (C_2 \times C_2) \cup (C_3 \times C_3) \cup (C_3 \times C_6) \cup (C_6 \times C_6)$,}
 \\
 - \dfrac{ 1 - x_i x_j }{ x_j - x_i }
 &\text{if $(i,j) \in (C_1 \times C_4) \cup (C_4 \times C_1) \cup (C_4 \times C_5)$,}
 \\
 \dfrac{ 1 - x_i x_j }{ x_j - x_i }
 &\text{if $(i,j) \in (C_2 \times C_3) \cup (C_2 \times C_6) \cup (C_3 \times C_2)$,}
 \\
 -1
 &\text{if $(i,j) \in \big( (C_1 \cup C_4) \times (C_2 \cup C_3 \cup C_6) \big) \cup (C_5 \times C_6)$,}
 \\
 1
 &\text{if $(i,j) \in \big( (C_2 \cup C_3) \times (C_1 \cup C_4 \cup C_5) \big) \cup (C_6 \times C_5)$.}
\end{cases}
$$
Let $\sigma$ be the ring automorphism of the Laurent polynomial ring 
$\Rat[x_1^{\pm 1}, \dots, x_n^{\pm 1}]$ defined by
$$
\sigma(x_i)
 = 
\begin{cases}
 x_i &\text{if $i \in I$,} \\
 x_i^{-1} &\text{if $i \not\in I$.}
\end{cases}
$$
Then, by using
\begin{equation}
\label{eq:sigma}
\sigma \left( \frac{ x_i - x_j }{ 1 - x_i x_j } \right)
 =
\begin{cases}
 \dfrac{ x_i - x_j }{ 1 - x_i x_j } &\text{if $(i,j) \in D^+_n(I) = I \times I \cup I^c \times I^c$,} \\
 \dfrac{ 1 - x_i x_j }{ x_i - x_j } &\text{if $(i,j) \in D^-_n(I) = I \times I^c \cup I^c \times I$,}
\end{cases}
\end{equation}
we have
$$
\sigma \left( X'_{i,j} \right)
 =
\begin{cases}
 - \dfrac{ x_j - x_i }{ 1 - x_i x_j }
 &\text{if $(i,j) \in (C_1 \cup C_4 \cup C_5) \times (C_1 \cup C_4 \cup C_5)$,}
\\
 \dfrac{ x_j - x_i }{ 1 - x_i x_j }
 &\text{if $(i,j) \in (C_2 \cup C_3 \cup C_6) \times (C_2 \cup C_3 \cup C_6)$,}
\\
 -1
 &\text{if $(i,j) \in (C_1 \cup C_4 \cup C_5) \times (C_2 \cup C_3 \cup C_6)$,}
\\
 1
 &\text{if $(i,j) \in (C_2 \cup C_3 \cup C_6) \times (C_1 \cup C_4 \cup C_6)$.}
\end{cases}
$$
Hence $\sigma \big( \Pf X'(I,J) \big) 
= \Pf \big( \sigma ( X'_{i,j} ) \big)_{1 \le i, j \le n}$ 
can be evaluated by the following lemma.

\begin{lemma}
\label{lem:PfY}
(\cite[Lemma~4.5]{Okada1998})
For a subset $K \subset [n]$, let $Y(K)$ be the $n \times n$ skew-symmetric matrix with $(i,j)$ entry 
given by
$$
Y(K)_{i,j}
 =
\begin{cases}
 - \dfrac{x_j-x_i}{1-x_ix_j} &\text{if $i \in K$ and $j \in K$,} \\
 \dfrac{x_j-x_i}{1-x_ix_j} &\text{if $i \not\in K$ and $j \not\in K$,} \\
 - 1 &\text{if $i \in K$ and $j \not\in K$,} \\
 1 &\text{if $i \not\in K$ and $j \in K$.} \\
\end{cases}
$$
Then we have
$$
\Pf Y(K)
 =
(-1)^{\Sigma(K)}
\prod_{(i,j) \in D^+_n(K)}
 \frac{x_j - x_i}{1 - x_i x_j},
$$
where $\Sigma(K) = \sum_{i \in K} i$ and
$$
D^+_n(K)
 =
\{ (i,j) \in [n] \times [n] : i < j \}
\cap
\big[ \big( K \times K \big) \cup \big( ([n] \setminus K) \times ([n] \setminus K) \big) \big].
$$
\end{lemma}

Applying this lemma to $K = C_1 \cup C_4 \cup C_5$, we have
$$
\Pf X'(I,J)
 =
(-1)^{\Sigma(C_1 \cup C_4 \cup C_5)}
\sigma \left(
 \prod_{(i,j) \in D^+_n(C_1 \cup C_4 \cup C_5)} \frac{x_j - x_i}{1 - x_i x_j}
\right).
$$
Then by using (\ref{eq:sigma}), we have
\begin{align}
\label{eq:L2}
\Pf X'(I,J)
 &=
(-1)^{\Sigma(C_1 \cup C_4 \cup C_5)}
\prod_{(i,j) \in T^{(1)}_{n,k} \cap \big( (C_1 \times C_1) \cup (C_2 \times C_2) \cup (C_3 \times C_3) \cup (C_4 \times C_4) \big)}
 \frac{x_j - x_i}{1 - x_i x_j}
\notag
\\
&\quad\times
\prod_{(i,j) \in T^{(1)}_{n,k} \cap \big( (C_1 \times C_4) \cup (C_2 \times C_3) \cup (C_3 \times C_2) \cup (C_4 \times C_1) \big)}
 \frac{1 - x_i x_j}{x_j - x_i}
\notag
\\
 &\quad\times
\prod_{(i,j) \in (C_1 \times C_5) \cup (C_3 \times C_6)}
 \frac{x_j - x_i}{1 - x_i x_j}
\prod_{(i,j) \in (C_2 \times C_6) \cup (C_4 \times C_5)}
 \frac{1 - x_i x_j}{x_j - x_i}
\notag
\\
&\quad\times
\prod_{(i,j) \in T^{(3)}_{n,k} \cap \big( (C_5 \times C_5) \cup (C_6 \cup C_6) \big)}
 \frac{x_j - x_i}{1 - x_i x_j}
\end{align}
Since $J = C_1 \cup C_3$, $[k] \setminus J = C_2 \cup C_4$ and
$$
f(x_i^{-1})
 = 
\frac{ (-1)^{n-k} }
     { \prod_{j \in C_5 \cup C_6} (x_j - x_i) },
\quad
f(x_i)
 =
\frac{ x_i^{n-k} }
     { \prod_{j \in C_5 \cup C_6} (1 - x_i x_j) },
$$
we see that
\begin{align*}
\prod_{i \in C_1} f(x_i^{-1})
\prod_{(i,j) \in C_1 \times C_5} \frac{x_j - x_i}{1 - x_i x_j}
 &=
(-1)^{(n-k) \# C_1}
\prod_{(i,j) \in C_1 \times C_5} \frac{ 1 }{1 - x_i x_j}
\prod_{(i,j) \in C_1 \times C_6} \frac{ 1 }{x_j - x_i},
\\
\prod_{i \in C_3} f(x_i^{-1})
\prod_{(i,j) \in C_3 \times C_6} \frac{x_j - x_i}{1 - x_i x_j}
 &=
(-1)^{(n-k) \# C_3}
\prod_{(i,j) \in C_3 \times C_6} \frac{ 1 }{1 - x_i x_j}
\prod_{(i,j) \in C_3 \times C_5} \frac{ 1 }{x_j - x_i},
\\
\prod_{i \in C_2} f(x_i)
\prod_{(i,j) \in C_2 \times C_6} \frac{1 - x_i x_j}{x_j - x_i}
&=
\prod_{i \in C_2} x_i^{n-k}
\prod_{(i,j) \in C_2 \times C_5} \frac{ 1 }{ x_j - x_i }
\prod_{(i,j) \in C_2 \times C_6} \frac{ 1 }{ 1 - x_i x_j},
\\
\prod_{i \in C_4} f(x_i)
\prod_{(i,j) \in C_4 \times C_5} \frac{1 - x_i x_j}{x_j - x_i}
&=
\prod_{i \in C_4} x_i^{n-k}
\prod_{(i,j) \in C_4 \times C_6} \frac{ 1 }{ x_j - x_i }
\prod_{(i,j) \in C_4 \times C_5} \frac{ 1 }{ 1 - x_i x_j}.
\end{align*}
Hence, combining these relations with (\ref{eq:L1}) and (\ref{eq:L2}), 
we obtain the following expression for the coefficient $L(I,J)$ in the left hand side:
\begin{align}
\label{eq:L}
L(I,J)
 &=
(-1)^{n-k + \Sigma(C_1 \cup C_4 \cup C_5) + (n-k) \# J}
\prod_{i \in [k] \setminus J} x_i^{n-k}
\notag
\\
 &\quad\times
\prod_{(i,j) \in T^{(1)}_{n,k} \cap \big( (C_1 \times C_1) \cup (C_2 \times C_2) \cup (C_3 \times C_3) \cup (C_4 \times C_4) \big)}
 \frac{x_j - x_i}{1 - x_i x_j}
\notag
\\
 &\quad\times
\prod_{(i,j) \in T^{(1)}_{n,k} \cap \big( (C_1 \times C_4) \cup (C_2 \times C_3) \cup (C_3 \times C_2) \cup (C_4 \times C_1) \big)}
 \frac{1 - x_i x_j}{x_j - x_i}
\notag
\\
 &\quad\times
\prod_{(i,j) \in \big( (C_1 \cup C_2) \times C_5 \big) \cup \big( (C_3 \cup C_4) \times C_6 \big)}
 \frac{ 1 }{1 - x_i x_j}
\prod_{(i,j) \in \big( (C_1 \cup C_2) \times C_6 \big) \cup \big( (C_3 \cup C_4) \times C_5 \big)}
 \frac{ 1 }{x_j - x_i}
\notag
\\
 &\quad\times
\prod_{(i,j) \in T^{(3)}_{n,k} \cap \big( (C_5 \times C_5) \cup (C_6 \times C_6) \big)}
 \frac{x_j - x_i}{1 - x_i x_j}.
\end{align}

Now we can finish the proof of Theorem~\ref{thm:Pf=det*det}.

\begin{demo}{Proof of Theorem~\ref{thm:Pf=det*det}}
We compare (\ref{eq:R}) with (\ref{eq:L}).
By the definitions (\ref{eq:D}) and (\ref{eq:C}), we have
\begin{align*}
T^{(1)}_{n,k} \cap D^+_n(I) \cap D^+_k(J)
 &=
T^{(1)}_{n,k} \cap \big( (C_1 \times C_1) \cup (C_2 \times C_2) \cup (C_3 \times C_3) \cup (C_4 \times C_4) \big),
\\
T^{(1)}_{n,k} \cap D^-_n(I) \cap D^-_k(J)
 &=
T^{(1)}_{n,k} \cap \big( (C_1 \times C_4) \cup (C_2 \times C_3) \cup (C_3 \times C_3) \cup (C_4 \times C_1) \big),
\\
T^{(2)}_{n,k} \cap D^+(I)
 &=
\big( (C_1 \cup C_2) \times C_5 \big) \cup \big( (C_3 \cup C_4) \times C_6 \big),
\\
T^{(2)}_{n,k} \cap D^-(I)
 &=
\big( (C_1 \cup C_2) \times C_6 \big) \cup \big( (C_3 \cup C_4) \times C_5 \big),
\\
T^{(3)}_{n,k} \cap D^+(I)
 &=
T^{(3)}_{n,k} \cap \big( (C_5 \times C_5) \cup (C_6 \times C_6) \big),
\end{align*}
and it follows that $R(I,J)$ coincides with $L(I,J)$ except for the sign.
Hence it remains to show that
\begin{equation}
\label{eq:sign}
k(k-1)/2 + n \# I + \Sigma(I) + k \# J + \Sigma(J)
 \equiv
n-k + \Sigma(C_1 \cup C_4 \cup C_5) + (n-k) \# J
\bmod 2.
\end{equation}
Since $I = C_1 \sqcup C_2 \sqcup C_5$, $J = C_1 \sqcup C_5$ and $[k] = C_1 \sqcup C_2 \sqcup C_3 \sqcup C_4$, 
we have
$$
\Sigma(I) + \Sigma(J) + \Sigma(C_1 \cup C_4 \cup C_5)
 =
2 \Sigma(C_1) + 2 \Sigma(C_5) + \Sigma([k])
 \equiv
k(k+1)/2
\bmod 2.
$$
Since $n$ is even, we obtain (\ref{eq:sign}).
Therefore we have $R(I,J) = L(I,J)$.
This completes the proof of Theorem~\ref{thm:Pf=det*det}, 
and hence of Theorem~\ref{thm:main}.
\end{demo}

\section{%
Shifted plane partitions of double staircase shape
}

In this section we find generating functions of shifted plane partitions 
by specializing the variables in the character identities in Theorem~\ref{thm:main}.
Also we derive Hopkins--Lai's formula for the number of lozenge tilings of flashlight regions.

\subsection{%
Generating functions of shifted plane partitions
}

For a strict partition $\mu$, 
we denote by $\AP(S(\mu))$ the set of all shifted plane partitions of shape $\mu$.
Given a shifted plane partition $\sigma \in \AP(S(\mu))$, 
we define the \emph{profile} $\pr(\sigma)$ to be the partition 
$(\sigma_{1,1}, \sigma_{2,2}, \dots)$ obtained by reading the main diagonal of $\sigma$.
For a set $\mathcal{Q}$ of partitions, we define 
$$
\AP(S(\mu);\mathcal{Q}) = \{ \sigma \in \AP(S(\mu)) : \pr(\sigma) \in \mathcal{Q} \}.
$$
For a subset $\mathcal{Q} \subset \Par((m^n))$, we write
$$
(a^n) + \mathcal{Q}
 =
\{ (a^n) + \lambda : \lambda \in \mathcal{Q} \}.
$$
Then, by specializing $x_i = q^i$ or $q^{i-1/2}$ in the character identities 
in Theorem~\ref{thm:main}, 
we obtain the generating functions for shifted plane partitions of shifted double staircase shape 
with respect to the weights $v(\sigma)$ and $w(\sigma)$ 
defined by (\ref{eq:wt1}) and (\ref{eq:wt2}) respectively.

\begin{theorem}
\label{thm:GF}
Suppose $0 \le k \le n$, and let $a$ and $m$ be nonnegative integers.
\begin{enumerate}
\item[(1)]
The generating functions for shifted plane partitions $\sigma$ of shape $\delta_n+\delta_k$ 
such that $(a^n) \subset \pr(\sigma) \subset ((a+m)^n)$ are given by
\begin{align}
&
\sum_{\sigma \in \AP(S(\delta_n+\delta_k) ; (a^n) + \Par((m^n)))} q^{v(\sigma)}
\notag
\\
 &\quad
=
q^{a n^2/2 - (m+2a) k^2/2}
\prod_{i=1}^n \frac{ [m/2+i-1/2] }{ [i-1/2] }
\prod_{1 \le i < j \le n} \frac{ [m+i+j-1] }{ [i+j-1] }
\notag
\\
&\quad\quad
\times
\prod_{i=1}^k \frac{ [m/2+a+i] }{ [i] }
\prod_{1 \le i < j \le k} \frac{ [m+2a+i+j] }{ [i+j] },
\label{eq:GF1a}
\\
&
\sum_{\sigma \in \AP(S(\delta_n+\delta_k) ; (a^n) + \Par((m^n)))} q^{w(\sigma)}
\notag
\\
 &\quad
=
q^{a n(n+1)/2 - (m+2a) k(k+1)/2}
\prod_{1 \le i \le j \le n} \frac{ [m+i+j-1] }{ [i+j-1] }
\prod_{1 \le i \le j \le k} \frac{ [m+2a+i+j] }{ [i+j] }.
\label{eq:GF1b}
\end{align}
\item[(2)]
If $m$ is even, then 
the generating functions for shifted plane partitions $\sigma$ of shape $\delta_n+\delta_k$ 
such that $\pr(\sigma) \in (a^n) + \Even((m^n))$ are given by
\begin{align}
&
\sum_{\sigma \in \AP(S(\delta_n+\delta_k) ; (a^n) + \Even((m^n)))} q^{v(\sigma)}
\notag
\\
 &\quad
=
q^{a n^2/2 - (m+2a) k^2/2}
\prod_{i=1}^n \frac{ [m/2+i] }{ [i] }
\prod_{1 \le i < j \le n} \frac{ [m+i+j] }{ [i+j] }
\notag
\\
 &\quad\quad
\times
\prod_{i=1}^k \frac{ [m/2+a+i] }{ [i] }
\prod_{1 \le i < j \le k} \frac{ [m+2a+i+j] }{ [i+j] },
\label{eq:GF2a}
\\
&
\sum_{\sigma \in \AP(S(\delta_n+\delta_k) ; (a^n) + \Even((m^n)))} q^{w(\sigma)}
\notag
\\
 &\quad
=
q^{a n(n+1)/2 - (m+2a) k(k+1)/2}
\prod_{1 \le i \le j \le n} \frac{ [m+i+j] }{ [i+j] }
\prod_{1 \le i \le j \le k} \frac{ [m+2a+i+j] }{ [i+j] }.
\label{eq:GF2b}
\end{align}
\item[(3)]
If $m > 0$, then 
the generating functions for shifted plane partitions $\sigma$ of shape $\delta_n+\delta_k$ 
such that $\pr(\sigma) \in (a^n) + \Even'((m^n))$ are given by
\begin{align}
&
\sum_{\sigma \in \AP(S(\delta_n+\delta_k) ; (a^n) + \Even'((m^n)))} q^{v(\sigma)}
\notag
\\
 &\quad
=
q^{a n^2/2 - (m+2a) k^2/2}
\prod_{i=1}^k \frac{ 1 }{ [2i-1] }
\prod_{1 \le i < j \le n} \frac{ [m+i+j-2] }{ [i+j-1] }
\prod_{1 \le i < j \le k} \frac{ [m+2a+i+j] }{ [i+j-1] }
\notag
\\
 &\quad\quad
\times
\left\{
 \prod_{i=1}^n \langle m/2+i-1 \rangle
 \prod_{i=1}^k [m/2+a+i]
\right.
\notag
\\
 &\quad\quad\quad\quad
+ (-1)^n
\left.
 (q-1)^{n-k}
 \prod_{i=1}^n [m/2+i-1]
 \prod_{i=1}^k \langle m/2+a+i \rangle
\right\},
\label{eq:GF3a}
\\
&
\sum_{\sigma \in \AP(S(\delta_n+\delta_k) ; (a^n) + \Even'((m^n)))} q^{w(\sigma)}
\notag
\\
 &\quad
=
q^{a n(n+1)/2 - (m+2a) k(k+1)/2}
\prod_{i=1}^k \frac{ 1 }{ [2i] }
\prod_{1 \le i < j \le n} \frac{ [m+i+j-2] }{ [i+j] }
\prod_{1 \le i < j \le k} \frac{ [m+2a+i+j] }{ [i+j] }
\notag
\\
 &\quad\quad
\times
\left\{
 \left( \sum_{l=0}^n q^{ml+l(l-1)} \qbinom{n}{l}^2 \right)
 \prod_{i=1}^k [m+2a+2i]
\right.
\notag
\\
 &\quad\quad\quad\quad
+ (-1)^n
\left.
 (q-1)^{n-k} \prod_{i=1}^n [m+2i-2]
 \left( \sum_{l=0}^k q^{(m+2a)l +l(l-1)} \qbinom{k}{l}^2 \right)
\right\},
\label{eq:GF3b}
\end{align}
where $\langle r \rangle = 1 + q^r$ and 
$$
\qbinom{n}{r}
 =
\frac{ [n] [n-1] \cdots [n-r+1] }
     { [r] [r-1] \cdots [1] }
$$
is the $q$-binomial coefficient.
\item[(4)]
If $m > 0$, then 
the generating functions for shifted plane partitions $\sigma$ of shape $\delta_n+\delta_k$ 
such that $\pr(\sigma) \in (a^n) + \Odd'((m^n))$ are given by
\begin{align}
&
\sum_{\sigma \in \AP(S(\delta_n+\delta_k) ; (a^n) + \Odd'((m^n)))} q^{v(\sigma)}
\notag
\\
 &\quad
=
q^{a n^2/2 - (m+2a) k^2/2}
\prod_{i=1}^k \frac{ 1 }{ [2i-1] }
\prod_{1 \le i < j \le n} \frac{ [m+i+j-2] }{ [i+j-1] }
\prod_{1 \le i < j \le k} \frac{ [m+2a+i+j] }{ [i+j-1] }
\notag
\\
 &\quad\quad
\times
\left\{
 \prod_{i=1}^n \langle m/2+i-1 \rangle
 \prod_{i=1}^k [m/2+a+i]
\right.
\notag
\\
 &\quad\quad\quad\quad
- (-1)^n
\left.
 (q-1)^{n-k}
 \prod_{i=1}^n [m/2+i-1]
 \prod_{i=1}^k \langle m/2+a+i \rangle
\right\},
\label{eq:GF4a}
\\
&
\sum_{\sigma \in \AP(S(\delta_n+\delta_k) ; (a^n) + \Odd'((m^n)))} q^{w(\sigma)}
\notag
\\
 &\quad
=
q^{a n(n+1)/2 - (m+2a) k(k+1)/2}
\prod_{i=1}^k \frac{ 1 }{ [2i] }
\prod_{1 \le i < j \le n} \frac{ [m+i+j-2] }{ [i+j] }
\prod_{1 \le i < j \le k} \frac{ [m+2a+i+j] }{ [i+j] }
\notag
\\
 &\quad\quad
\times
\left\{
 \left( \sum_{l=0}^n q^{ml+l(l-1)} \qbinom{n}{l}^2 \right)
 \prod_{i=1}^k [m+2a+2i]
\right.
\notag
\\
 &\quad\quad\quad\quad
- (-1)^n
\left.
 (q-1)^{n-k} \prod_{i=1}^n [m+2i-2]
 \left( \sum_{l=0}^k q^{(m+2a)l +l(l-1)} \qbinom{k}{l}^2 \right)
\right\}.
\label{eq:GF4b}
\end{align}
\end{enumerate}
\end{theorem}

Some special cases of this theorem appeared in the earlier literatures.
By specializing $q=1$, Equations (\ref{eq:GF1a}) and (\ref{eq:GF1b}) with $a=0$ 
reduce to Hopkins--Lai's product formula \cite[Theorem~1]{HL2021} (see Theorem~\ref{thm:HL}) 
for the number of shifted plane partitions of shifted double staircase shape.
Putting $k=0$ and $a=0$ in (\ref{eq:GF1a}) and (\ref{eq:GF1b}), 
we obtain the MacMahon and Bender--Knuth (ex-)conjectures 
(\ref{eq:MacMahon}) and (\ref{eq:BenderKnuth}) mentioned in Introduction respectively.
And we can recover the formulas of \cite[Theorem~1, cases (CYH) and (CYI)]{Proctor1990} 
by specializing $k=0$ and $a=0$ in (\ref{eq:GF2a}) and (\ref{eq:GF2b}).
The special case $p=n$ of \cite[Theorems~11 and 21]{Krattenthaler1995} 
are obtained by considering the case $k=0$ and $a=1$ of (\ref{eq:GF2a}) and (\ref{eq:GF2b}).
Similarly we derive the special cases ($p=0$ and $p=c$) of \cite[Theorem~6]{Krattenthaler1998} 
from (\ref{eq:GF3a})--(\ref{eq:GF4b}) with $k=0$ and $a=0$.
By considering the case $k=n-1$ or $k=n$, we obtain the generating functions 
for shifted plane partitions of trapezoidal shape.
For example, (\ref{eq:GF1a}) and (\ref{eq:GF1b}) are $q$-analogues of 
(special cases of) the formula of \cite[Theorem~1]{Proctor1983}.
And we can recover special cases of \cite[Theorems~4 and 7]{Krattenthaler1998} 
from (\ref{eq:GF2a})--(\ref{eq:GF4b}) with $k=n-1$ or $n$.

\subsection{%
Proof of Theorem~\ref{thm:GF}
}

In this subsection we derive Theorem~\ref{thm:GF} from Theorem~\ref{thm:main}.

The following bijection enables us to convert the $(k,n-k)$-symplectic characters 
into the generating functions for shifted plane partitions of shape $\delta_n+\delta_k$. 
(This bijection is a generalization of that used in \cite[Lemma~2]{Proctor1983} 
and \cite[Proof of Theorem~4]{Krattenthaler1998}.)

\begin{lemma}
\label{lem:bijection}
For a partition $\lambda$ of length $\le n$, 
let $\AP(S(\delta_n+\delta_k);\lambda)$ be the set of shifted plane partitions 
of shifted double staircase shape $\delta_n+\delta_k$ with profile $\lambda$.
Then there is a bijection $\phi : \AP(S(\delta_n+\delta_k);\lambda) 
\to \Tab^{(k,n-k)}(\lambda)$ satisfying
\begin{equation}
\label{eq:wt}
\vectx^{\phi(\sigma)} 
\Big|_{\substack{x_i=q^{k-i+1/2} \, (1 \le i \le k) \\ x_i=q^{n+k-i+1/2} \, (k+1 \le i \le n)}}
 = q^{v(\sigma)},
\quad
\vectx^{\phi(\sigma)} 
\Big|_{\substack{x_i=q^{k-i+1} \, (1 \le i \le k) \\ x_i=q^{n+k-i+1} \, (k+1 \le i \le n)}}
 = q^{w(\sigma)},
\end{equation}
for $\sigma \in \AP(S(\delta_n+\delta_k);\lambda)$, 
where the weights $v(\sigma)$ and $w(\sigma)$ is defined in (\ref{eq:wt1}) 
and (\ref{eq:wt2}) respectively.
In particular, we have
\begin{gather}
\label{eq:GF=symp-a}
\sum_{\sigma \in \AP(S(\delta_n+\delta_k);\lambda)} q^{v(\sigma)}
 =
\symp^{(k,n-k)}_\lambda (q^{1/2}, q^{3/2}, \dots, q^{n-1/2}),
\\
\label{eq:GF=symp-b}
\sum_{\sigma \in \AP(S(\delta_n+\delta_k);\lambda)} q^{w(\sigma)}
 =
\symp^{(k,n-k)}_\lambda (q, q^2, \dots, q^n).
\end{gather}
\end{lemma}

\begin{demo}{Proof}
To a shifted plane partition $\sigma \in \AP(S(\delta_n+\delta_k);\lambda)$, 
we associate a column-strict plane partition $\pi$ of shape $\lambda$ 
whose $i$th row is the conjugate partition of the $i$th row of $\sigma$.
Then, by replacing $1, 2, \dots, n-k-1, n-k,\dots, n+k$ with 
$n, n-1, \dots, k+1, \overline{k}, k, \dots, \overline{1}, 1$ respectively, 
we obtain a $(k,n-k)$-symplectic tableau $T \in \Tab^{(k,n-k)}(\lambda)$.
For $n=4$, $k=2$ and $\lambda = (4,3,1,1)$, 
an example is given by
$$
\sigma
 =
\begin{matrix}
 4 & 4 & 2 & 2 & 1 & 0 \\
   & 3 & 2 & 2 & 1 \\
   &   & 1 & 1 \\
   &   &   & 1
\end{matrix}
\longmapsto
\pi =
\begin{matrix}
 5 & 4 & 2 & 2 \\
 4 & 3 & 1 \\
 2 \\
 1
\end{matrix}
\longmapsto
T =
\begin{matrix}
 \overline{1} & 2 & 3 & 3 \\
 2 & \overline{2} & 4 \\
 3 \\
 4
\end{matrix}.
$$
Since these procedures are invertible, the correspondence $\sigma \mapsto T$ gives 
a bijection between $\AP(S(\delta_n+\delta_k);\lambda)$ and $\Tab^{(k,n-k)}(\lambda)$.
And, since the multiplicity of $l$ in $\pi$ is equal to 
the difference of the traces $t_{l-1}(\sigma) - t_l(\sigma)$, 
where $t_l(\sigma)$ is defined by (\ref{eq:tr}) and $t_{n+k}(\sigma) = 0$, 
the multiplicity $m_T(\gamma)$ of $\gamma$ in $T$ is given by
\begin{align*}
m_T(i)
 &=
t_{n+k-2i+1}(\sigma) - t_{n+k-2i+2}(\sigma)
\quad(1 \le i \le k),
\\
m_T(\overline{i})
 &=
t_{n+k-2i}(\sigma) - t_{n+k-2i+1}(\sigma)
\quad(1 \le i \le k),
\\
m_T(i)
 &=
t_{n-i}(\sigma) - t_{n-i+1}(\sigma)
\quad(k+1 \le i \le n).
\end{align*}
Then we can show
\begin{align*}
v(\sigma)
 &=
\sum_{i=1}^k \left( k-i+\frac{1}{2} \right) \left( m_T(i) - m_T(\overline{i}) \right)
 +
\sum_{i=k+1}^n \left( n+k-i+\frac{1}{2} \right) m_T(i),
\\
w(\sigma)
 &=
\sum_{i=1}^k (k-i+1) \left( m_T(i) - m_T(\overline{i}) \right)
 +
\sum_{i=k+1}^n (n+k-i+1) m_T(i),
\end{align*}
which imply (\ref{eq:wt}).
It follows from the Jacobi--Trudi-type identity (\ref{eq:JT2}) that
$$
\symp^{(k,n-k)}_\lambda(x_1, \dots, x_k | x_{k+1}, \dots, x_n)
 =
\symp^{(k,n-k)}_\lambda(x_k, \dots, x_1 | x_n, \dots, x_{k+1}).
$$
By using this symmetry, we can obtain (\ref{eq:GF=symp-a}) and (\ref{eq:GF=symp-b}).
\end{demo}

In order to derive Theorem~\ref{thm:GF}, 
we need the following formulas for the specializations 
of classical group characters corresponding to rectangular partitions.

\begin{lemma}
\label{lem:specialization}
\begin{enumerate}
\item[(1)]
If we specialize $x_i = q^{i-1/2}$ for $1 \le i \le n$, then we have
\begin{gather}
\label{eq:symp_qa}
\symp_{(m^n)}(q^{1/2}, q^{3/2}, \dots, q^{n-1/2})
 =
\frac{ 1 }{ q^{m n^2/2} }
\prod_{i=1}^n \frac{ [m+i] }{ [i] }
\prod_{1 \le i < j \le n} \frac{ [2m+i+j] }{ [i+j] },
\\
\label{eq:oddorth_qa}
\orth^B_{(m^n)}(q^{1/2}, q^{3/2}, \dots, q^{n-1/2})
 =
\frac{ 1 }{ q^{m n^2/2} }
\prod_{i=1}^n \frac{ [m+i-1/2] }{ [i-1/2] }
\prod_{1 \le i < j \le n} \frac{ [2m+i+j-1] }{ [i+j-1] },
\\
\label{eq:evenorth_qa}
\orth^D_{(m^n)}(q^{1/2}, q^{3/2}, \dots, q^{n-1/2})
 =
\frac{ \chi(m) }{ q^{m n^2/2} }
\prod_{i=1}^n \frac{ \langle m+i-1 \rangle }{ \langle i-1 \rangle },
\prod_{1 \le i < j \le n} \frac{ [2m+i+j-2] }{ [i+j-2] },
\end{gather}
where $\chi(m) = 1$ if $m=0$ and $2$ if $m>0$.
\item[(2)]
If we specialize $x_i = q^i$ for $1 \le i \le n$, then we have
\begin{gather}
\label{eq:symp_qb}
\symp_{(m^n)}(q, q^2, \dots, q^n)
 =
\frac{ 1 }{ q^{m n(n+1)/2} }
\prod_{1 \le i \le j \le n} \frac{ [2m+i+j] }{ [i+j] },
\\
\label{eq:oddorth_qb}
\orth^B_{(m^n)}(q, q^2, \dots, q^n)
 =
\frac{ 1 }{ q^{m n(n+1)/2} }
\prod_{1 \le i \le j \le n} \frac{ [2m+i+j-1] }{[i+j-1] }.
\\
\label{eq:evenorth_qb}
\orth^D_{(m^n)}(q, q^2, \dots, q^n)
 =
\frac{ 1 }{ q^{m n(n+1)/2} }
\prod_{1 \le i < j \le n} \frac{ [2m+i+j-2] }{ [i+j] }
\sum_{l=0}^n q^{2ml+l(l-1)} \qbinom{n}{l}^2.
\end{gather}
\end{enumerate}
\end{lemma}

\begin{demo}{Proof}
As the proofs are similar, we only give a proof of (\ref{eq:evenorth_qb}).
It follows from the definition (\ref{eq:evenorth_bialt}) that
$$
\orth^D_{(m^n)}(q, q^2, \dots, q^n)
 =
\chi(m)
\frac{ \det \Big( (q^{m+j-1})^i + (q^{m+j-1})^{-i} \Big)_{1 \le i, j \le n} }
     { \det \Big( (q^{j-1})^i + (q^{j-1})^{-i} \Big)_{1 \le i, j \le n} }.
$$
Here we use the following determinant evaluation (see \cite[(7.6)]{Krattenthaler1998}):
$$
\det \left( z_i^j + z_i^{-j} \right)_{1 \le i, j \le n}
 =
(z_1 \dots z_n)^{-n}
\prod_{1 \le i < j \le n} (z_i - z_j)(1 - z_i z_j)
\sum_{l=0}^n e_l(z_1, \dots, z_n)^2.
$$
Then we have
\begin{multline*}
\det \Big( (q^{m+j-1})^i + (q^{m+j-1})^{-i} \Big)_{1 \le i, j \le n}
\\
 =
\prod_{i=1}^n (q^{m+i-1})^i
\prod_{1 \le i < j \le n} ( 1 - q^{j-i} ) ( 1 - q^{2m+i+j-2} )
\sum_{l=0}^n \left( q^{lm + l(l-1)/2} \qbinom{n}{l} \right)^2.
\end{multline*}
Hence by using
$$
\sum_{l=0}^n e_l(1, q, q^2, \dots, q^{n-1})^2
 =
q^{n(n-1)/2}
e_n(1, q, \dots, q^{n-1}, 1, q^{-1}, \dots, q^{-(n-1)})
 =
2 \prod_{i=1}^n \frac{ [n+i-1] }{ [i] },
$$
we can complete the proof of (\ref{eq:evenorth_qb}).
\end{demo}

Now we are ready to prove Theorem~\ref{thm:GF}.

\begin{demo}{Proof of Theorem~\ref{thm:GF}}
Identities (\ref{eq:GF1a})--(\ref{eq:GF2b}) are immediately obtained by specializing 
$x_i = q^{i-1/2}$ or $q^i$ in Theorem~\ref{thm:main} and using Lemmas~\ref{lem:bijection} 
and \ref{lem:specialization}.
We need a few more manipulations to obtain (\ref{eq:GF3a})--(\ref{eq:GF4b}).
Here we give a proof of (\ref{eq:GF3a}).
By using the character identities (\ref{eq:main3}) and (\ref{eq:main4}), we obtain
\begin{align*}
&
\sum_{\lambda \in \Even'((m^n))} \symp^{(k,n-k)}_{(a^n)+\lambda}(\vectx)
\\
 &=
\frac{1}{2}
\Bigg(
 \orth^D_{((m/2)^n)}(x_1, \dots, x_n)
 \cdot
 \symp_{((m/2+a)^k)}(x_1, \dots, x_k)
 \cdot
 (x_{k+1} \dots x_n)^{m/2+a}
\\
 &\quad\quad\quad
 +
 (-1)^n
 \symp_{((m/2-1)^n)}(x_1, \dots, x_n)
 \cdot
 \orth^D_{((m/2+a+1)^k)}(x_1, \dots, x_k)
 \cdot
 \prod_{i=k+1}^n 
 x_i^{m/2+a} (x_i-x_i^{-1})
\Bigg).
\end{align*}
Specializing $x_i = q^{i-1/2}$ ($1 \le i \le n$) and using
$$
\prod_{1 \le i < j \le r} \frac{ 1 }{ [i+j-2] }
\prod_{i=1}^r \frac{ 1 }{ \langle i-1 \rangle }
 =
\prod_{1 \le i < j \le r} \frac{ 1 }{ [i+j-1] },
$$
we can arrive at the expression (\ref{eq:GF3a}).
\end{demo}

\subsection{%
Lozenge tilings of flashlight regions
}

In this subsection, we apply the character identity (Theorem~\ref{thm:main} (1)) to 
enumerate the lozenge tilings of flashlight regions in the triangular lattice.

\begin{figure}[htb]
\centering
\includegraphics[bb=0 0 243 207]{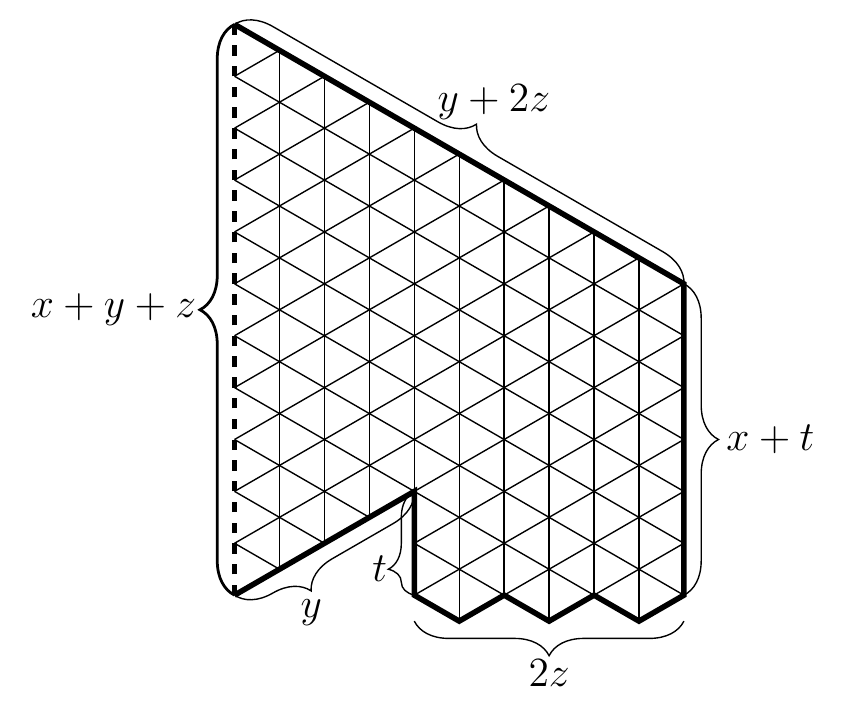}
\caption{The flashlight region $F_{x,y,z,t}$ for $x=4$, $y=4$, $z=3$, $t=2$}
\label{fig:2}
\end{figure}

\begin{figure}[htb]
\centering
\includegraphics[bb=0 0 175 52]{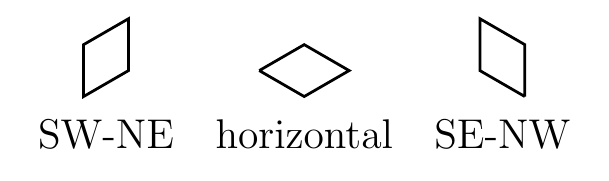}
\caption{Three types of lozenges}
\label{fig:3}
\end{figure}

A lozenge is the union of two unit equilateral triangles joined along an edge.
A lozenge tiling of a region $R$ in the regular triangular lattice 
is a covering of $R$ by lozenges with neither gap nor overlap.
For nonnegative integers $x$, $y$, $z$ and $t$, 
let $F_{x,y,z,t}$ be the ``flashlight region'' shown in Figure~\ref{fig:2}, 
where the dashed line indicates a free boundary, 
i.e., lozenges are allowed to protrude across it.
We consider lozenge tilings of $F_{x,y,z,t}$ by three types of lozenges given in Figure~\ref{fig:3}.
Then Hopkins and Lai \cite{HL2021} obtain the following product formula for the number of 
lozenge tilings of $F_{x,y,z,t}$, 
and use the case $t=0$ to derive the formula for the number of shifted plane partitions 
of shifted double staircase shape (Theorem~\ref{thm:HL}).

\begin{theorem}
\label{thm:tiling}
(\cite[Theorem~1.2]{HL2021} for $y>0$)
For nonnegative integers $x$, $y$, $z$ and $t$, 
the number $M(F_{x,y,z,t})$ of lozenge tilings of the region $F_{x,y,z,t}$ is given by
\begin{equation}
\label{eq:tiling}
M(F_{x,y,z,t})
 =
\prod_{1 \le i \le j \le y+z}
 \frac{ x+i+j+1 }
      { i+j+1 }
\prod_{1 \le i \le j \le z}
 \frac{ x+2t+i+j }
      { i+j }.
\end{equation}
\end{theorem}

The proof by Hopkins--Lai \cite{HL2021} is based on an extension of Kuo condensation 
to regions with a free boundary due to Ciucu \cite{Ciucu2020}.
Here we use the character identity in Theorem~\ref{thm:main} (1) to give an alternate proof.

\begin{figure}[htb]
\centering
\includegraphics[bb=0 0 188 150]{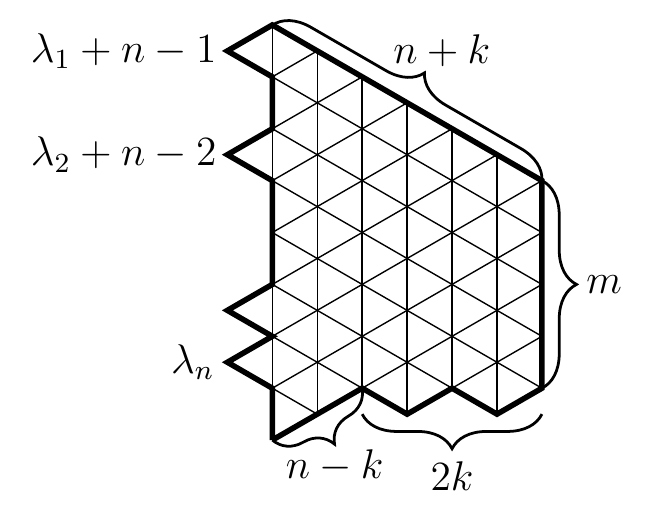}
\caption{The region $R^{(k,n-k)}_m(\lambda)$ for $k=2$, $n-k=2$, $m=4$, $\lambda = (4,3,1,1)$}
\label{fig:4}
\end{figure}

\begin{figure}[htb]
$$
\raisebox{-60pt}{
\includegraphics[bb=0 0 104 163]{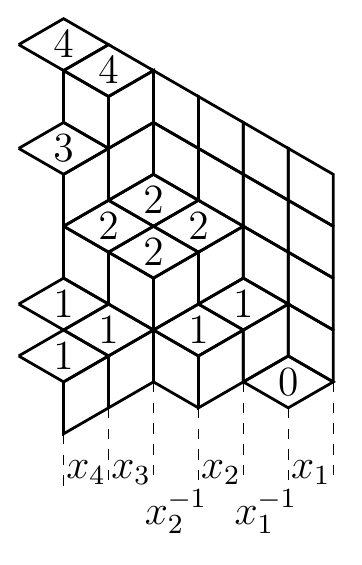}
}
\longleftrightarrow
\begin{matrix}
 4 & 4 & 2 & 2 & 1 & 0 \\
   & 3 & 2 & 2 & 1 \\
   &   & 1 & 1 \\
   &   &   & 1
\end{matrix}
\longleftrightarrow
\begin{matrix}
 \overline{1} & 2 & 3 & 3 \\
 2 & \overline{2} & 4 \\
 3 \\
 4
\end{matrix}
$$
\caption{The correspondence among tilings, shifted plane partitions and intermediate symplectic tableuax}
\label{fig:5}
\end{figure}

In order to apply the character identity, 
we give an interpretation of intermediate symplectic characters 
in terms of lozenge tilings.
Let $0 \le k \le n$.
Given a partition $\lambda$ such that $l(\lambda) \le n$ and $\lambda_1 \le m$,
we denote by $R^{(k,n-k)}_m(\lambda)$ the region obtained from $F_{m,n-k,k,0}$ 
by adjoining $n$ left-pointing triangles at the positions 
$\lambda_1+n-1, \dots, \lambda_{n-1}+1, \lambda_n$ 
to the free boundary of $F_{m,n-k,k,0}$, 
where the vertical edges on the free boundary are labeled with $0, 1, \dots, m+n-1$ 
from bottom to top.
For example, the region $R^{(2,2)}_4 (4,3,1,1)$ is depicted in Figure~\ref{fig:4}.
We denote by $\TT(R^{(k,n-k)}_m(\lambda))$ the set of lozenge tilings of $R^{(k,n-k)}(\lambda)$.
For a tiling $T$ of $R^{(k,n-k)}_m(\lambda)$, or $F_{m,n-k,k,a}$, 
we define its weight $\wt(T)$ as follows.
We assign to each SW-NE lozenge $L$ in the $i$th column from the right a weight
$$
\wt(L)
 =
\begin{cases}
x_j &\text{if $1 \le i \le 2k$ and $i=2j-1$,} \\
x_j^{-1} &\text{if $1 \le i \le 2k$ and $i=2j$,} \\
x_{i-k} &\text{if $2k+1 \le i \le n+k$,}
\end{cases}
$$
and put $\wt(L) = 1$ for the other two types of lozenges $L$.
Then the weight $\wt(T)$ of a tiling $T$ is defined 
as the product of all the weights of the lozenges used in the tiling $T$. 
For example, if $T$ is the tiling given in the left picture of Figure~\ref{fig:5}, 
then its weight is computed as
$$
\wt (T)
 = (x_1)^0 \cdot (x_1^{-1})^1 \cdot (x_2)^2 \cdot (x_2^{-1})^1 \cdot (x_3)^3 \cdot (x_4)^2
 = x_1^{-1} x_2 x_3^3 x_4^2.
$$
The following lemma provides an interpretatioin of $\symp^{(k,n-k)}_\lambda(\vectx)$ 
as the tiling generating function, 
and interpolates between \cite[Theorem~2.3]{AF2020} (the Schur function case, $k=0$) 
and \cite[Theorem~2.8]{AF2020} (the symplectic character case, $k=n$), 
both of which are given in terms of perfect matchings.

\begin{lemma}
\label{lem:sp-tiling}
For a partitin $\lambda$ of length $\le n$ with $\lambda_1 \le m$, 
there is a weight-preserving bijection 
between $\Tab^{(k,n-k)}(\lambda)$ and $\TT(R^{(k,n-k)}_m(\lambda))$.
Hence we have 
$$
\symp^{(k,n-k)}_\lambda(x_1, \dots, x_k|x_{k+1}, \dots, x_n)
 =
\sum_{T \in \TT(R^{(k,n-k)}_m(\lambda))} \wt(T).
$$
\end{lemma}

\begin{demo}{Proof}
By Lemma~\ref{lem:bijection}, we have a bijection between $\Tab^{(k,n-k)}(\lambda)$ 
and $\AP^m(S(\delta_n+\delta_k);\lambda)$. 
And there is a natural bijection between 
$\TT(R^{(k,n-k)}_m(\lambda))$ and $\AP^m(S(\delta_n+\delta_k);\lambda)$, 
which is obtained by reading the ``height'' of the horizontal lozenges 
along the paths consisting of horizontal and SW-NE lozenges.
The desired bijection is obtained by composing these two bijections.
See Figure~\ref{fig:5}.
\end{demo}

Now Theorem~\ref{thm:tiling} follows from Theorem~\ref{thm:main} (1) 
by using Lemma~\ref{lem:sp-tiling}.

\begin{demo}{Proof of Theorem~\ref{thm:tiling}}
Let $M_{m,n-k,k,a} = M_{m,n-k,k,a}(x_1, \dots, x_k | x_{k+1}, \dots, x_n)$ 
be the generating function of tilings of $F_{m,n-k,k,a}$.
We prove 
\begin{equation}
\label{eq:tiling-gf}
M_{m,n-k,k,a}
 =
\orth^B_{(m/2)^n}(x_1, \dots, x_n)
\cdot
\symp_{(m/2+a)^k}(x_1, \dots, x_k)
\cdot
(x_{k+1} \cdots x_n)^{m/2}.
\end{equation}
By specializing $x_1 = \dots = x_n = 1$ in (\ref{eq:tiling-gf}) 
and using (\ref{eq:symp_qa}), (\ref{eq:oddorth_qa}) with $q=1$, 
we obtain the desired identity (\ref{eq:tiling}).

Let $\tilde{F}_{m,n-k,k,a}$ be the region obtained from $F_{m+a,n-k,k,0}$ by changing 
$a$ edges from the bottom on the free boundary into a non-free boundary (see Figure~\ref{fig:6}), 
\begin{figure}[htb]
\centering
\includegraphics[bb=0 0 238 215]{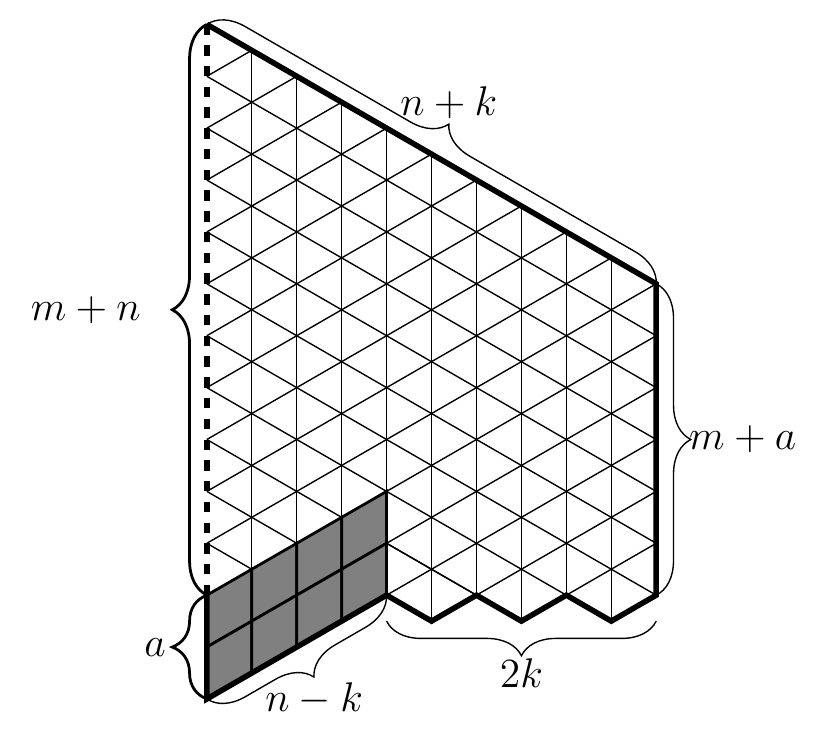}
\caption{The region $\tilde{F}_{m,n-k,k,a}$ for $m=4$, $n-k=4$, $k=3$, $a=2$}
\label{fig:6}
\end{figure}
and $\tilde{M}_{m,n-k,k,a}$ the tiling generating function of $\tilde{F}_{m,n-k,k,a}$.
Because of the presence of the non-free boundary, 
the bottom-left corner (the shaded region in Figure~\ref{fig:6}) is tiled 
by forced SW-NE lozenges, and we have
$$
\tilde{M}_{m,n-k,k,a}
 =
(x_{k+1} \cdots x_n)^a
\cdot
M_{m,n-k,k,a}.
$$
On the other hande, it follows from Lemma~\ref{lem:sp-tiling} that 
$$
\tilde{M}_{m,n-k,k,a}
 =
\sum_{\lambda \in (a^n) + \Par((m^n))} 
 \symp^{(k,n-k)}_\lambda(x_1, \dots, x_k|x_{k+1}, \dots, x_n).
$$
Hence Equation~(\ref{eq:tiling-gf}) follows from (\ref{eq:main1}) in Theorem~\ref{thm:main}.
\end{demo}

\begin{remark}
The above proof of Theorem~\ref{thm:tiling} 
shows that the character identity (\ref{eq:main1}) is equivalent 
to the tiling generating function identity (\ref{eq:tiling-gf}).
When $k<n$, one can prove (\ref{eq:tiling-gf}) by a method similar to that of the proof of 
\cite[Theorem~1.2]{HL2021}, 
i.e., by using the induction on $m+k$ 
and appealing to Ciucu's generalization \cite[Corollary~1]{Ciucu2020} of Kuo condensation.
\end{remark}
\appendix
\section{%
More determinant formulas for $\symp^{(k,n-k)}_\lambda$
}

In this appendix, we apply the theory of Macdonald's ninth variation of Schur functions 
\cite{Macdonald1992} 
to derive dual Jacobi--Trudi and Giambelli formulas 
for intermediate symplectic characters.

Recall Macdonald's ninth variation of Schur functions.
Let $\{ h^{[p]}_r : p, r \in \Int, r \ge 1 \}$ be indeterminates.
We use the convention $h^{[p]}_0 = 1$ and $h^{[p]}_r = 0$ for $r < 0$.
For a partition $\lambda$ we define
\begin{equation}
\label{eq:9th_schur}
s^{[p]}_\lambda = \det \Big( h^{[p+j-1]}_{\lambda_i-i+j} \Big)_{1 \le i, j \le l(\lambda)}.
\end{equation}
And we write $e^{[p]}_r = s^{[p]}_{(1^r)}$ for $r \ge 0$ 
and put $e^{[p]}_r = 0$ for $r < 0$.
Then we have

\begin{prop}
\label{prop:9th_schur}
\begin{enumerate}
\item[(1)]
(\cite[(9.6')]{Macdonald1992})
If $\lambda' = (\lambda'_1, \dots, \lambda'_m)$ is the conjugate partition of $\lambda$, 
then we have
\begin{equation}
\label{eq:9th_dualJT}
s^{[p]}_\lambda
 =
\det \Big(
 e^{[p-j+1]}_{\lambda'_i-i+j}
\Big)_{1 \le i , j \le m}.
\end{equation}
\item[(2)]
(\cite[(9.7)]{Macdonald1992})
If $\lambda = (\alpha_1, \dots, \alpha_r|\beta_1, \dots, \beta_r)$ in the Frobenius notation, 
then we have
\begin{equation}
\label{eq:9th_G}
s^{[p]}_\lambda
 = 
\det 
\Big(
 s^{[p]}_{(\alpha_i|\beta_j)}
\Big)_{1 \le i, j \le r}.
\end{equation}
\end{enumerate}
\end{prop}

Now we specialize
$$
h^{[p]}_r
 =
\begin{cases}
 h_r(x_p^{\pm 1}, \dots, x_k^{\pm 1}, x_{k+1}, \dots, x_n)
 &\text{if $p \le k$,} \\
 h_r(x_p, \dots, x_n)
 &\text{if $k+1 \le p \le n$,} \\
 0
 &\text{if $p \ge n+1$,}
\end{cases}
$$
where $(\dots, x_{-2}, x_{-1}, x_0, x_1, \dots, x_k, x_{k+1}, \dots, x_n)$ are indeterminates.
Then, by comparing (\ref{eq:JT1}) and (\ref{eq:9th_schur}), we have
$$
s^{[p]}_\lambda
 =
\begin{cases}
 \symp^{(k-p+1,n-k)}_\lambda(x_p, \dots, x_k|x_{k+1}, \dots, x_n)
 &\text{if $l(\lambda) \le n-p+1$ and $p \le k$,} \\
 s_\lambda(x_p, \dots, x_n)
 &\text{if $l(\lambda) \le n-p+1$ and $k+1 \le p \le n$,} \\
 0 &\text{otherwise.}
\end{cases}
$$
Hence the following Giambelli formula for intermediate symplectic characters 
is an immediate consequence of (\ref{eq:9th_G}).
(See \cite[p.29]{FK1997} for a combinatorial proof in the case $l(\lambda) \le k+1$.)

\begin{prop}
\label{prop:G}
If a partition $\lambda$ of length $\le n$ is written as 
$\lambda = (\alpha_1, \dots, \alpha_r|\beta_1, \dots, \beta_r)$ in the Frobenius notation, 
then we have
\begin{equation}
\label{eq:G}
\symp^{(k,n-k)}_\lambda(x_1, \dots, x_k|x_{k+1}, \dots, x_n)
 =
\det \left(
 \symp^{(k,n-k)}_{(\alpha_i|\beta_j)} (x_1 \dots, x_k|x_{k+1}, \dots, x_n)
\right)_{1 \le i, j \le r}.
\end{equation}
\end{prop}

To state dual Jacobi--Trudi formulas, we introduce some notations.
We define $e^\circ_r(x_1^{\pm 1}, \dots, x_k^{\pm 1})$ by the generating function
\begin{equation}
\label{eq:gen_e0}
\sum_{r \ge 0} e^\circ_r(x_1^{\pm 1}, \dots, x_k^{\pm 1}) t^r
 =
(1 - t^2) \prod_{i=1}^k (1 + x_i t) (1 + x_i^{-1} t).
\end{equation}
Then 
$e^\circ_r(x_1^{\pm 1}, \dots, x_k^{\pm 1}) = \symp_{(1^r)}(x_1, \dots, x_k)$ for $0 \le r \le k$, and
$$
\symp^{(k,n-k)}_{(1^r)}(x_1, \dots, x_k|x_{k+1}, \dots, x_n)
 =
\sum_{p=0}^k
 e^\circ_p(x_1^{\pm 1}, \dots, x_k^{\pm 1}) e_{r-p}(x_{k+1}, \dots, x_n).
$$
For an integer $m$, we define
$$
e^{(k,n-k)}_{r,m}(x_1, \dots, x_k | x_{k+1}, \dots, x_n)
 =
\sum_{p=0}^{k-m}
 e^\circ_p(x_1^{\pm 1}, \dots, x_k^{\pm 1}) e_{r-p}(x_{k+1}, \dots, x_n).
$$
With these notations we have the following dual Jacobi--Trudi formulas.

\begin{prop}
\label{prop:dualJT}
Let $\lambda$ be a partition of length $\le n$ and $\lambda'$ the conjugate partition of $\lambda$.
\begin{enumerate}
\item[(1)]
We have
\begin{multline}
\label{eq:dualJT1}
\symp^{(k,n-k)}_\lambda(x_1, \dots, x_k | x_{k+1}, \dots, x_n)
\\
 =
\det \left(
 \sum_{m=0}^{j-1} e^{(k,n-k)}_{\lambda'_i-i+j-2m,m}(x_1, \dots, x_k|x_{k+1}, \dots, x_n)
\right)_{1 \le i, j \le \lambda_1}.
\end{multline}
\item[(2)]
(\cite[Proposition~8.1]{Proctor1991})
If $l(\lambda) \le k+1$, then we have
\begin{multline}
\label{eq:dualJT2}
\symp^{(k,n-k)}_\lambda(x_1, \dots, x_k | x_{k+1}, \dots, x_n)
\\
 =
\det \left(
 \begin{cases}
  e^{(k,n-k)}_{\lambda'_i-i+1,0}(x_1, \dots, x_k | x_{k+1}, \dots, x_n)
  &\text{if $j=1$} \\
  e^{(k,n-k)}_{\lambda'_i-i+j,0}(x_1, \dots, x_k | x_{k+1}, \dots, x_n) \\
  \quad + e^{(k,n-k)}_{\lambda'_i-i-j+2,0}(x_1, \dots, x_k | x_{k+1}, \dots, x_n)
  &\text{if $2 \le j \le \lambda_1$}
 \end{cases}
\right)_{1 \le i, j \le \lambda_1}.
\end{multline}
\end{enumerate}
\end{prop}

If $k=0$, then $e^{(0,n)}_{r,m} = 0$ for $m \ge 1$ and $e^{(0,n)}_{r,0} = e_r$, 
hence (\ref{eq:dualJT1}) reduces to 
the dual Jacobi--Trudi formula for Schur functions:
$$
s_\lambda(x_1, \dots, x_n)
 =
\det \left( e_{\lambda'_i-i+j}(x_1, \dots, x_n) \right)_{1 \le i, j \le \lambda_1}.
$$
And, if $k=n$, then (\ref{eq:dualJT2}) gives the dual Jacobi--Trudi formula for 
symplectic characters:
$$
\symp_\lambda(x_1, \dots, x_n)
 =
\det \left(
 \begin{cases}
  e^\circ_{\lambda'_i-i+1}(x_1^{\pm 1}, \dots, x_n^{\pm 1})
  &\text{if $j=1$} \\
  e^\circ_{\lambda'_i-i+j}(x_1^{\pm 1}, \dots, x_n^{\pm 1})
  + e^\circ_{\lambda'_i-i-j+2}(x_1^{\pm 1}, \dots, x_n^{\pm 1})
  &\text{if $2 \le j \le \lambda_1$}
 \end{cases}
\right)_{1 \le i, j \le \lambda_1}.
$$

In the proof of Proposition~\ref{prop:dualJT}, 
we use the following relations:

\begin{lemma}
\label{lem:rel_e}
Let $\vecty = (x_1, \dots, x_k)$, $\vectz = (x_{k+1}, \dots, x_n)$ and $u$ be indeterminates.
Then we have
\begin{equation}
\label{eq:rel_e1}
e^{(k+1,n-k)}_{r,m}(u,\vecty|\vectz)
 =
e^{(k,n-k)}_{r,m-1}(\vecty|\vectz)
 + e^{(k,n-k)}_{r-2,m+1}(\vecty|\vectz)
 + (u+u^{-1}) e^{(k,n-k)}_{r-1,m}(\vecty|\vectz),
\end{equation}
and $e^{(k,n-k)}_{r,-1}(\vecty|\vectz) = e^{(k,n-k)}_{r,0}(\vecty|\vectz)$.
If $ r \le k+1$, then we have
\begin{equation}
\label{eq:rel_e2}
e^{(k+1,n-k)}_{r,0}(u,\vecty|\vectz)
 =
e^{(k,n-k)}_{r,0}(\vecty|\vectz)
 + e^{(k,n-k)}_{r-2,0}(\vecty|\vectz)
 + (u+u^{-1}) e^{(k,n-k)}_{r-1,0}(\vecty|\vectz).
\end{equation}
\end{lemma}

\begin{demo}{Proof}
By the generating function (\ref{eq:gen_e0}), we have
$$
e^\circ_r (u^{\pm 1}, x_1^{\pm 1}, \dots, x_k^{\pm 1})
 =
e^\circ_r (x_1^{\pm 1}, \dots, x_k^{\pm 1})
 + e^\circ_{r-2} (x_1^{\pm 1}, \dots, x_k^{\pm 1})
 + (u+u^{-1}) e^\circ_{r-1} (x_1^{\pm 1}, \dots, x_k^{\pm 1}),
$$
and $e^\circ_{k+1}(x_1^{\pm 1}, \dots, x_k^{\pm 1}) = 0$.
The claims follows easily from these relations.
\end{demo}

\begin{demo}{Proof of Proposition~\ref{prop:dualJT}}
By using (\ref{eq:9th_dualJT}), we have
\begin{equation}
\label{eq:dualJT0}
\symp^{(k,n-k)}_\lambda(\vecty|\vectz)
 =
\det \left(
 e^{(k+j-1,n-k)}_{\lambda'_i-i+j,0}(x_{2-j}, \dots, x_0, x_1, \dots, x_k | x_{k+1}, \dots, x_n)
\right)_{1 \le i, j \le \lambda_1},
\end{equation}
where $(x_{2-n}, \dots, x_{-2}, x_{-1}, x_0)$ are dummy indeterminates.
In order to prove the proposition, we perform column operations 
on the matrix in the right hand side of (\ref{eq:dualJT0}).

(1)
Let $\vecty = (y_1, \dots, y_k)$ and $\vectz = (z_1, \dots, z_{n-k})$ be indeterminates.
Given an integer sequence $\alpha = (\alpha_1, \dots, \alpha_n)$ and a nonnegative integer $r$, 
let $C^{(k,n-k)}_{\alpha,r}(\vecty|\vectz)$ be the column vector with the $i$th entry
$$
\sum_{m=0}^r
 e^{(k,n-k)}_{\alpha_i+r-2m,m}(\vecty|\vectz).
$$
Consider the $n \times (l+2)$ matrix $G$ given by
$$
G
 = 
\begin{pmatrix}
C^{(k,n-k)}_{\alpha,0}(\vecty|\vectz)
 &
C^{(k+1,n-k)}_{\alpha+1,0}(u,\vecty|\vectz)
 &
C^{(k+1,n-k)}_{\alpha+1,1}(u,\vecty|\vectz)
 &
\cdots
 &
C^{(k+1,n-k)}_{\alpha+1,l}(u,\vecty|\vectz)
\end{pmatrix},
$$
where $\alpha+1 = (\alpha_1+1, \dots, \alpha_n+1)$ and $u$ is another indeterminate.
Now we perform the following column operations:
\begin{enumerate}
\item[(i)]
subtract the $1$st column multiplied by $u+u^{-1}$ from the $2$nd column;
\item[(ii)] 
subtract the $2$nd column multiplied by $u+u^{-1}$ from the $3$rd column, 
and then subtract the $1$st column from the $3$rd column;
\item[(iii)]
subtract the $3$nd column multiplied by $u+u^{-1}$ from the $4$rd column, 
and then subtract the $2$nd column from the $4$rd column;
\item[(iv)]
and so on.
\end{enumerate}
Then, by using the relation (\ref{eq:rel_e1}), the matrix $G$ is transformed into
$$
\tilde{G}
 =
\begin{pmatrix}
C^{(k,n-k)}_{\alpha,0}(\vecty|\vectz)
 &
C^{(k,n-k)}_{\alpha,1}(\vecty|\vectz)
 &
C^{(k,n-k)}_{\alpha,2}(\vecty|\vectz)
 &
\cdots
 &
C^{(k,n-k)}_{\alpha,l+1}(\vecty|\vectz)
\end{pmatrix},
$$
where we note that the additional variable $u$ is eliminated.
By repeated application of this process, 
we can transform the matrix on the right hand side of (\ref{eq:dualJT0}) 
into the matrix on the right hand side of (\ref{eq:dualJT1}).
This completes the proof of (\ref{eq:dualJT1}).

(2) can be proved in a similar fashion by using the relation (\ref{eq:rel_e2}).
\end{demo}



\begin{thebibliography}{99}

\bibitem{Andrews1977a}
G.~E.~Andrews,
MacMahon's conjecture on symmetric plane partitions,
Proc. Nat. Acad. Sci. U.S.A. {\bfseries 74} (1977), 426--429.

\bibitem{Andrews1978}
G.~E.~Andrews,
Plane partitions. I. The MacMahon conjecture,
Adv. in Math. Suppl. Stud. {\bfseries 1} (1978), 131--150.

\bibitem{Andrews1977b}
G.~E.~Andrews,
Plane partitions. II. The equivalence of the Bender-Knuth and MacMahon conjectures, 
Pacific J. Math. {\bfseries 72} (1977), 283--291.

\bibitem{AB2019}
A.~Ayyer and R.~E.~Behrend,
Factorization theorems for classical group characters, 
with applications to alternating sign matrices and plane partitions,
J. Combin. Theory Ser. A {\bfseries 165} (2019), 78---105.

\bibitem{AF2020}
A.~Ayyer and I.~Fischer,
Bijective proofs of skew Schur polynomial factorizations,
J. Combin. Theory Ser. A {\bfseries 174} (2020), 105241.

\bibitem{BK1972}
E.~A.~Bender and D.~E.~Knuth,
Enumeration of plane partitions.
J. Combin. Theory Ser. A {\bfseries 13} (1972), 40--54.

\bibitem{Ciucu2020}
M.~Ciucu,
Correlation of a macroscopic dent in a wedge with mixed boundary conditions, 
Trans. Amer. Math. Soc. {\bfseries 373} (2020), 2173--2190.

\bibitem{FK1997}
M.~Fulmek and C.~Krattenthaler,
Lattice path proofs for determinant formulas for symplectic and orthogonal characters,
J. Combin. Theory Ser. A {\bfseries 77} (1997), 3--50.

\bibitem{FH}
W.~Fulton and J.~Harris,
``Representation theory, A first course'',
Grad. Texts in Math. {\bfseries 129}, Springer-Verlag, New York, 1991.

\bibitem{GV1985}
I.~Gessel and G.~Viennot,
Binomial determinants, paths, and hook length formulae,
Adv. Math. {\bfseries 58} (1985), 300--321.

\bibitem{GV1989}
I.~Gessel and G.~Viennot,
Determinants, paths, and plane partitions,
preprint.

\bibitem{Gordon1983}
B.~Gordon, 
A proof of the Bender-Knuth conjecture,
Pacific J. Math. {\bfseries 108} (1983), 99--113.

\bibitem{HL2021}
S.~Hopkins and T.~Lai,
Plane partitions of shifted double staircase shape,
J. Combin. Theory Ser. A {\bfseries 183} (2021), 105486.

\bibitem{IW1995}
M.~Ishikawa and M.~Wakayama,
Minor summation formula of Pfaffians,
Linear and Multilinear Algebra, {\bfseries 39} (1995), 285--305.

\bibitem{King1976}
R.~C.~King,
Weight multiplicities for the classical groups,
in ``Group Theoretical Methods in Physics (Fourth Internat. Colloq., Nijmegen, 1975)'', 
Lecture Notes in Phys. {\bfseries 50}, Springer, Berlin, 1976, 
pp. 490--499. 

\bibitem{KT1987}
K.~Koike and I.~Terada,
Young-diagrammatic methods for the representation theory of the classical groups 
of type $B_n$, $C_n$, $D_n$, 
J. Algebra {\bfseries 107} (1987), 466--511.

\bibitem{Krattenthaler1995}
C.~Krattenthaler,
The major counting of nonintersecting lattice paths and generating functions for tableaux,
Mem. Amer. Math. Soc. {\bfseries 115} (1995), no. 552, vi+109 pp.

\bibitem{Krattenthaler1998}
C.~Krattenthaler,
Identities for classical group characters of nearly rectangular shape,
J. Algebra {\bfseries 209} (1998), 1--64.

\bibitem{Lindstrom1973}
B.~Lindstr\"om,
On the vector representations of induced matroids,
Bull. London Math. Soc. {\bfseries 5} (1973), 85--90.

\bibitem{Macdonald1992}
I.~G.~Macdonald,
Schur functions: Theme and variations,
S\'em. Lothar. Combin. {\bfseries 28} (1992), 5--39.

\bibitem{Macdonald1995}
I.~G.~Macdonald,
``Symmetric Functions and Hall Polynomials, 2nd edition'',
Oxford Univ. Press, 1995.

\bibitem{MacMahon1899}
P.~A.~MacMahon,
Partitions of numbers whose graphs possess symmetry,
Trans. Cambridge Philos. Soc. {\bfseries 17} (1898), 149--170.

\bibitem{Okada1998}
S.~Okada,
Applications of minor summation formulas to rectangular-shaped
representations of classical groups,
J. Algebra {\bfseries 205} (1998), 337--367.

\bibitem{Okada2020}
S.~Okada,
A bialternant formula for odd symplectic characters and its application,
Josai Math. Monographs {\bfseries 12} (2020), 99--116.

\bibitem{Proctor1983}
R.~A.~Proctor,
Shifted plane partitions of trapezoidal shape, 
Proc. Amer. Math. Soc. {\bfseries 89} (1983), 553--559.

\bibitem{Proctor1988}
R.~A.~Proctor,
Odd symplectic groups, 
Invent. Math. {\bfseries 92} (1988), 307--332

\bibitem{Proctor1990}
R.~A.~Proctor,
New symmetric plane partition identities from invariant theory 
work of De Concini and Procesi,
European J. Combin. {\bfseries 11} (1990), 289--300.

\bibitem{Proctor1991}
R.~A.~Proctor,
A generalized Berele--Schensted algorithm 
and conjectured Young tableaux for intermediate symplectic groups,
Trans. Amer. Math. Soc. {\bfseries 324} (1991), 655--692.

\bibitem{Proctor1994}
R.~A.~Proctor,
Young tableaux, Gelfand patterns, and branching rules for classical groups, 
J. Algebra {\bfseries 164} (1994), 299--360.

\end{thebibliography}
\end{document}